\documentclass{article}
\usepackage[utf8]{inputenc}
\usepackage[T1]{fontenc}
\usepackage{lmodern}

\usepackage{amsmath}
\usepackage{amsthm}
\usepackage{amssymb}
\usepackage{mathtools}
\usepackage{mathrsfs}

\usepackage[ruled,vlined]{algorithm2e}

\usepackage{pgf,tikz,pgfplots}
\usetikzlibrary{arrows}
\usepackage{pstricks-add}
\usepackage{pst-func}
\usepackage{graphicx}
\usepackage{color}
\usepackage{xcolor}

\usepackage{enumitem}
\usepackage{multirow}
\usepackage{makecell}
\usepackage{float}

\usepackage[a4paper, margin=1in]{geometry}
\usepackage{hyperref}
\usepackage{cleveref}

\usepackage{authblk}

\theoremstyle{plain}
\newtheorem{theorem}{Theorem}[section]
\newtheorem{proposition}[theorem]{Proposition}
\newtheorem{lemma}[theorem]{Lemma}

\newtheorem{corollary}[theorem]{Corollary}

\newtheorem{remark}[theorem]{Remark}

\theoremstyle{definition}
\newtheorem{definition}[theorem]{Definition}

\everymath{\displaystyle}

\newcommand{\sign}{\operatorname{sign}}

\newcommand{\Matrix}[1]{\begin{bmatrix} #1 \end{bmatrix}}
\newcommand{\diag}[1]{\operatorname{diag}\circbrack{#1}}

\newcommand{\R}{\mathbb{R}}

\newcommand{\Z}{\mathbb{Z}}
\newcommand{\Q}{\mathbb{Q}}
\newcommand{\Tr}{\mathrm{Tr}}
\newcommand{\GL}{\mathrm{GL}}
\newcommand{\Nm}{\mathrm{N}}
\newcommand{\Gal}{\mathrm{Gal}}

\newcommand{\set}[1]{\left\{#1\right\}}

\newcommand{\abs}[1]{\left|#1\right|}

\newcommand{\circbrack}[1]{\left(#1\right)}
\newcommand{\T}{\mathfrak{t}}

\author[1]{Nam H. Le}
\author[2]{Dat T. Tran}
\author[3]{David Karpuk}
\author[4]{Ha T. N. Tran}

\affil[1]{Department of Mathematics and Statistics, Florida Atlantic University, US}
\affil[2]{Department of Mathematics, University of Science, Ho Chi Minh City, Vietnam }
\affil[3]{QMill, Espoo, Finland} %david.karpuk@qmill.com}
\affil[4]{Augustana Faculty, University of Alberta,  Canada}

\title{Well-Rounded Twists of the Ring of Integers in Cyclic Cubic Fields}
\newcommand{\keywords}[1]{\par\noindent\textbf{Keywords: } #1}
\newcommand{\subjclass}[1]{\par\noindent\textbf{MSC (2020): } #1}

\begin{document}
	\maketitle

\abstract{Computing well-rounded twists of ideals in number fields has been done when the field degree is $2$. In this paper, we develop a new algorithm to detect whether a basis of an ideal $\mathfrak{I}$ in a cyclic cubic field $F$ yields a well-rounded twist of $\mathfrak{I}$. We then prove that under certain conditions on a given basis of the ring of integers $\mathcal{O}_F$, the existence of its well-rounded twist is equivalent to the existence of a principal well-rounded ideal in $K$. Applying the result and the algorithm, we explicitly compute well-rounded twists of the ring of integers for cyclic cubic fields in the families of Shanks, Washington, and Kishi. In addition, we show that infinitely many fields in Shanks's family have rings of integers that admit orthogonal well-rounded twists. }
    
	%%%%%%%%%%%%%%%%%%%%%%%%%%%%%%
\vspace{0.2cm}


\subjclass{Primary:
        11Y40, %Algebraic number theory computations
        11R16, %Cubic and quartic extensions
        11R04 %(1980–now)Algebraic numbers; rings of algebraic integers
        6B10, %Lattice ideals, congruence relations; 	
	Secondary: 
	11H06. %Lattices and convex bodies (number-theoretic aspects)
	% https://mathscinet.ams.org/msnhtml/msc2020.pdf
   }	

    \section{Introduction}

    \subsection{Background and Related Work}
    
    A \emph{lattice} is a discrete additive subgroup of Euclidean space. Lattices play a significant role in many areas of mathematics, with extensive connections to discrete optimization and Lie group theory, and are central in the theory of sphere-packings \cite{conway2013sphere}.  Lattice-based codebooks are used in communication over channels with additive noise \cite[Chapter III]{conway2013sphere}, and more recently, lattices arising from totally real number fields have found applications in wireless communication \cite{oggier2004algebraic,fading1996good,vehkalahti2013inverse}. 

    A lattice of full-rank in Euclidean space is called \textit{well-rounded} if its set of shortest vectors spans the entire space. The concept of well-roundedness goes back to at least \cite{ash1984}, but has seen recent interest due to applications of such lattices in communications \cite{damir2019well,gnilke2016well}.  
    
    The simplest example of a well-rounded lattice is the orthogonal lattice $\mathbb{Z}^n$.  In \cite{oggier2004algebraic}, the authors show how to apply a \emph{twist} map to a lattice $L$ arising from an ideal of the ring of integers of a totally real field to construct rotated versions of $\mathbb{Z}^n$ which have so-called \emph{full-diversity}, meaning that no non-zero lattice point has any coordinate equal to zero.  More generally, one can apply a twisting map to any lattice and hope that the result is well-rounded, an idea that was first explored in \cite{mcmullen2005minkowski,solan2019stable,levin2017closedorbits}, wherein the authors show various results establishing the existence of well-rounded twists of algebraic and arbitrary lattices.  Explicit methods for constructing all well-rounded twists of ideal lattices were given in \cite{damir2019well} for the case of real quadratic fields, and in \cite{le2022well} for imaginary quadratic fields. However, explicit methods for constructing well-rounded twists of lattices (even ideal lattices) in $\mathbb{R}^3$ are largely non-existent.  It is the goal of the current work to make progress on this question. 

    \subsection{Summary of Current Contribution}
    
    In the case of a full-rank $L$ lattice in $\R^2$, given a $\Z$-basis of $L$, one can easily determine whether this basis is minimal (see \Cref{def:minimal_basis}) or not by checking the angle between the two vectors in the basis (see \cite[Proposition 1.1]{fukshansky2012integral}). Such criteria also exist for lattices in $\R^3$ and are undoubtedly known to experts, but we make them explicit for the sake of making the current work self-contained.  To identify a well-rounded twist of a lattice $L\subseteq\R^3$, one now chooses a twist map $T$ and a basis $B$, and verifies minimality for the twisted basis $TB$.  If the twist map $T$ is such that $TL$ is well-rounded with minimal basis $TB$, we say that $B$ is a \emph{good} basis (see \Cref{def:good_basis}).  As every well-rounded lattice has a minimal basis, this reduces the problem to enumerating all good bases $B$, and goodness can be verified by checking a simple list of explicit algebraic equations.  \Cref{sec:prelim_result} of the current paper is devoted to the above line of argument. We develop a new algorithm \Cref{rem:twist_ideal_cyclic} to detect whether a basis of an ideal $\mathfrak{I}$ in a cyclic cubic field $F$ yields a well-rounded twist of $\mathfrak{I}$.
    
    % Secondly, a generalization of well-rounded lattices is well-rounded twist lattices (see \cite[Section 1]{solan2019stable} and \Cref{def:twist}). We call a basis $B$ of a lattice $L$ a \emph{good basis}  if there exists a twist map $T$ such that $TL$ is well-rounded with $TB$ a minimal basis. In general, and hence in the cyclic cubic case, every ideal lattice admits a well-rounded twist \cite[Theorem 1.2]{solan2019stable}, however, in the cubic case we can go further and given explicit conditions for a basis to be good.
    
    % Finally, in light of recent results on well-rounded lattices in a cyclic cubic field $F$ \cite{tran2023well}, a natural question arises: is there a connection between well-rounded twists of lattice given by the ring of integers $\mathcal{O}_F$, and the well-rounded ideal lattices presented in \cite[Section 3]{tran2023well}? \dk{I don't really understand this paragraph.}\HT{Dave, please feel free to delete/modify it}
    
% In this paper, we address these three questions for ideal lattices \HT{we only do WR twists for the ring of integers, right? } arising from cyclic cubic fields.

The remaining sections are dedicated to instantiating the above results when our lattice comes from the ring of integers $\mathcal{O}_F$ of a real cyclic cubic field $F$ via the canonical embedding.  Next, We show that if $F$ admits a well-rounded principal ideal $\mathfrak{I}$, then $\mathcal{O}_F$ has a well-rounded twist similar to $\mathfrak{I}$. For the converse, we show that under certain conditions, a good basis of $\mathcal{O}_F$ implies the existence of a principal well-rounded ideal whose norm divides the discriminant of $F$ (\Cref{prop:wash_good_bases}). We prove that every minimal basis of such a WR ideal $\mathfrak I$ must be of the form $\set{\pm \kappa, \pm \sigma\circbrack{\kappa}, \pm \sigma^2\circbrack{\kappa}}$, where $\kappa$ is a shortest vector of $\mathfrak I$ by Propositions \ref{prop:con_basis_3midm} and \ref{prop:con_basis_3nmidm}. In such cases, the Gram matrix of the twisted basis (\Cref{def:twist_twisted}) has off-diagonal entries equal in absolute value, giving a concrete criterion under which the converse holds. We also show that $\mathcal{O}_F$ admits an orthogonal well-rounded twist when its different ideal is principal (\Cref{thm:ortho_twist} and \Cref{prop:cyc_cu_orth}). In particular, there are infinitely many simplest cubic fields from Shanks \cite{shanks1974simplest} of which the ring of integers has an orthogonal WR twist lattice (see \Cref{prop:infinite_ortho}). 
We conclude the paper by providing explicit good bases for three special families of cyclic cubic fields from Shanks ~\cite{shanks1974simplest}, Washington \cite{washington1997family} and Kishi \cite{kishi2003family}, applying the results in Propositions \ref{prop:wr_non_ortho_shank}, \ref{prop:wash_good_bases}, and \ref{prop:kishi_good_bases}.

    % \dk{We basically outline the structure of the paper three different times.  Will work on condensing this.} \HT{yes, please}

%Our main contribution is 

%This paper is structured as follows. In \Cref{sec:prelim_result}, we begin by recalling the definitions of lattices and well-rounded lattices in $\R^3$. We then prove a necessary and sufficient condition for a lattice of rank 3 to be well-rounded. Next, we discuss twists of lattices in $\R^3$ and determine the form of twist matrices in this case. We also review ideal lattices, defining polynomials, integral bases, and discriminants of cyclic cubic fields, along with several established properties and lemmas that will be used in later sections. In \Cref{sec:ideal_norm_div_dis}, we examine the well-roundedness of all ideal lattices  whose norms divide the discriminant of a cyclic cubic field $F$. \Cref{sec:wr_twist_of_roi} presents the main results of the paper on WR twists of the ring of integers of cyclic cubic fields. In \Cref{sec:explicit_wr_twists}, we construct explicit good bases for several well-known families of cyclic cubic fields. We conclude with a summary of our findings and suggestions for future research in \Cref{sec:conclusion}.  \dk{I say we can delete this paragraph entirely but I'll leave it up to you.}

Additionally, to support some of the computational steps in the proofs, we used \texttt{SageMath}~\cite{sagemath} 
to compute minimal polynomials, discriminants, and Gram matrices, 
and \texttt{Mathematica}~\cite{Mathematica} to compute and simplify norms, 
as well as to solve inequalities arising in \Cref{thm:wr3_nec_suff} and \Cref{lem:sign_of_three_numbers}. 
All source code and computational outputs are available in  GitHub repository~\cite{github_computations_wr}.

	%%%%%%%%%%%%%%%%%%%%%%%%%%%%%%
	\section{Preliminary results}\label{sec:prelim_result}

	\subsection{Lattices and Well-Rounded Lattices in $\R^3$}\label{sec:wr_lattice}

	Let $\mathcal{B}=\{v_1,v_2,...,v_m\}$ be a linearly independent set of vectors in $\mathbb{R}^n$, $1\le m \le n$. The set $L=\left\{\sum_{i=1}^{m}a_iv_i|a_i\in \mathbb{Z}\right\}$ is called a \textit{lattice} in $\mathbb{R}^n$ of rank $m$ and the set $\mathcal{B}$ is said to be a \textit{basis} of $L$. In case $m=n$, we say that $L$ is a \textit{full-rank lattice}. In this paper, we only consider full-rank lattices.

Let $L_1, L_2$ be two full-rank lattices in $\R^n$. We say that $L_1$ and $L_2$ are \textit{similar} if there exists a nonzero real constant $c$ and an $n\times n$ real orthogonal matrix $O$ such that $L_2 = cOL_1$. 
    
	The value $|L| = \min_{0\ne u \in L}\|u\|^2$ is called the \textit{smallest norm} or \textit{the first minimum} of the lattice $L\subseteq \mathbb{R}^n$, where $\|.\|$ denotes the usual Euclidean norm in $\mathbb{R}^n$, and the set of \textit{shortest vectors} of $L$ is defined as 
	\[S(L):=\{u\in L:\|u\|^2=|L|\}.\]
	
	\begin{definition} \label{def:WR}
		A lattice $L$ in $\mathbb{R}^n$
		is called \textit{well-rounded (WR)} if $S(L)$ generates $\mathbb{R}^n$, that is, if $S(L)$ contains $n$ linearly independent vectors.
	\end{definition}
	\begin{definition}\label{def:minimal_basis}
		A basis $B$ of a WR lattice in $\R^n$ is called a \textit{minimal basis} if $B \subseteq S(L)$.
	\end{definition}
	\begin{proposition}\label{min_vector_basis}
		Let $\Lambda$ be a WR lattice of rank $k \leq 3$.  Then any $k$ linearly independent minimal vectors of $\Lambda$ form a basis of $\Lambda$.
	\end{proposition}
	\begin{proof}
		Let $\Lambda'\subseteq\Lambda$ be a sublattice spanned by any $k$ linearly independent minimal vectors of $\Lambda$.  Then by \cite[Corollary 2.6.10]{martinet2013perfect} we have $[L:L']\leq \gamma_k^{k/2}$, where $\gamma_k$ is the Hermite constant in dimension $k$.  Since $\gamma_k^{k/2} < 2$ for $k = 2,3$ we must have $[L:L'] = 1$ and the result follows.
	\end{proof}
	
	In this paper, we only consider lattices in $\R^3$. Let $\Lambda$ be a lattice of rank $3$.  To any basis $B = \{x_1, x_2, x_3\}$ of $\Lambda$ we can associate the Gram matrix $G_B = \left(\langle x_i, x_j\rangle \right)$ where $\langle x_i, x_j\rangle$ is the usual inner (scalar) product of $x_i$ and $x_j$.  When $n = 3$ we will usually denote a basis $B$ by $\{x,y,z\}$, and we will set
    \begin{equation}\label{eq:inner_prod}
        u = \langle x,y\rangle, \qquad v = \langle x,z\rangle, \qquad w = \langle y,z\rangle.
    \end{equation}
	
	The following result is undoubtedly known to experts, but we include an independent proof for the sake of the current work being self-contained.
	
	\begin{theorem}\label{thm:wr3_nec_suff}
		Let $\Lambda$ be a lattice of rank $3$.  Then $\Lambda$ is WR if and only if there exists a basis $B = \{x,y,z\}$ of $\Lambda$ whose associated Gram matrix is of the form
		\[
		G_B = \begin{bmatrix}
			s & u & v \\
			u & s & w \\
			v & w & s
		\end{bmatrix}
		\]
		where $s = \|x\|^2 = \|y\|^2 = \|z\|^2$ and the quantities $s,u,v,w$ defined as in \eqref{eq:inner_prod} satisfy the inequalities
		\begin{align*}
			\max\{ |u|, |v|, |w| \}&\leq s/2, \text{  and } \\
			\max\{    -u + v +w,  u - v +w, u + v - w,   -u - v - w\}&\leq s.
		\end{align*}
        Moreover, such a basis $B$ is minimal. Conversely, in case $\Lambda$ is WR, the Gram matrix of any minimal basis of $\Lambda$ must satisfy the above conditions.  
	\end{theorem}
	\begin{proof}
		First, suppose that $\Lambda$ is WR.  By \Cref{min_vector_basis} we can find a basis $B = \{x,y,z\}$ consisting of minimal vectors of $\Lambda$.  Since the sublattice of rank $2$ spanned by any pair of $x,y,z$ must also be WR, we must have $|u|,|v|,|w|\leq s/2$.  Now consider the set of eight vectors $\mathcal{S} = \{\pm x \pm y \pm z\}\subseteq\Lambda$.  We have $\|t\|^2 \leq s$ for all $t\in \mathcal{S}$.  Setting $t = x + y + z$ and expanding out $\|t\|^2 \leq s$ yield the inequality $-u-v-w\leq s$.  Negating the coefficient of any of $x,y,z$ in $t$ will negate exactly two of the quantities $u,v,w$, giving the remaining three desired inequalities.
		
		Now suppose that $\Lambda$ has a basis $B = \{x,y,z\}$ satisfying the stated inequalities.  Let $t = \alpha_1x + \alpha_2y + \alpha_3z\in\Lambda$ for some integers $\alpha_i$, and suppose that $\|t\|^2 < s$. We will show that $t = 0$, which implies that the vectors in $B$ are minimal and thus $\Lambda$ is WR.  This is equivalent to
		\[
		(\alpha_1^2 + \alpha_2^2 + \alpha_3^2)s + 2\alpha_1\alpha_2u + 2\alpha_1\alpha_3v + 2\alpha_2\alpha_3w < s.
		\]
		We will proceed by showing that $t=0$ on a case-by-case basis according to the signs of the $\alpha_i$.
		
		First, suppose that all $\alpha_i \geq 0$. Then the inequalities $|u|,|v|,|w|\leq s/2$ imply that
		\[
		(\alpha_1^2 + \alpha_2^2 + \alpha_3^2)s- \alpha_1\alpha_2s - \alpha_1\alpha_3s - \alpha_2\alpha_3s \leq (\alpha_1^2 + \alpha_2^2 + \alpha_3^2)s + 2\alpha_1\alpha_2u + 2\alpha_1\alpha_3v + 2\alpha_2\alpha_3w. 
		\]
		Combining these inequalities and dividing by the positive number $s$ yields
		\[
		(\alpha_1^2 + \alpha_2^2 + \alpha_3^2) - \alpha_1\alpha_2 - \alpha_1\alpha_3 - \alpha_2\alpha_3 < 1
		\]
		or equivalently,
		\[
		(\alpha_1 - \alpha_2)^2 + (\alpha_1 - \alpha_3)^2 + (\alpha_2 - \alpha_3)^2 < 2.
		\]
		The left-hand side of this inequality must be an even integer and hence is equal to zero, from which we conclude that $\alpha_1 = \alpha_2 = \alpha_3$.  Denoting this quantity by $\alpha$, we have $w = \alpha(x + y + z)$.  Using the inequality $- u - v -w \leq s$, we compute that
		\[
		s > \|t\|^2 = \alpha^2(3s + 2u + 2v + 2w) \geq \alpha^2 s.
		\]
		This forces $\alpha = 0$ and hence $t = 0$.  
		
		Now suppose that $\alpha_1 \geq 0$ and $\alpha_2,\alpha_3\leq 0$.  A similar calculation as in the previous case shows that the inequalities $|u|,|v|,|w|\leq s/2$ imply $\alpha_1 = -\alpha_2 = -\alpha_3$.  Thus we can write $w = \alpha(-x + y + z)$ for some integer $\alpha$.  The inequality $u + v - w \leq s$ now implies $s > \|w\|^2 = \alpha^2(3s-2u-2v+2w)\geq \alpha^2 s$ from which we again conclude that $\alpha = 0$ and thus $w = 0$.
		
		The remainder of the proof follows as above with two more cases of signs for the $\alpha_i$ to check, each corresponding to one of the as-of-yet unused inequalities in the statement of the theorem. We conclude that $\Lambda$ has a basis consisting of minimal vectors and is therefore WR.

        The last two statements of the theorem follow from the previous proof and  \Cref{def:minimal_basis}.\\	
	\end{proof}

	\subsection{Twists of lattices in $\R^3$}\label{sec:twist_lattices}
	Let us define the \emph{diagonal group} to be
	\[
	\mathcal{A}_3 = \left\{
	\begin{bmatrix}
		\alpha & 0 &  0\\
		0& \beta & 0\\
		0 & 0 & \gamma
	\end{bmatrix} :\ \alpha,\beta,\gamma \in \mathbb{R}_{>0}, \alpha \beta \gamma =1
	\right\}
	\subseteq \GL_3(\R).
	\]
        \begin{definition} \label{def:twist_twisted}   Let $\Lambda$ be a lattice in $\mathbb{R}^3$ with a basis $B$, and let $T \in \mathcal{A}_3$.
        \begin{itemize}
            \item The \textit{twist} of the lattice $\Lambda$ by $T$ is defined as the lattice $T\Lambda$, obtained by multiplying $T$ to each vector in $\Lambda$, that is
    $T\Lambda=\{Tv\mid v\in \Lambda\}.$
\item      We also denote by \(B\) the matrix whose columns are the basis vectors in \(B\). The lattice \(T\Lambda\) has a basis given by the columns of the matrix \(TB\). We call $TB$ the\textit{ twisted basis} of $B$ by $T$.
        \end{itemize}
        \end{definition}
   
        \begin{definition}\label{def:twist}
    Let $T \in \mathcal{A}_3$. A twist $T\Lambda$ is called a \textit{WR twist} of $\Lambda$ if the lattice  $T\Lambda$ is WR.
        \end{definition}
        \begin{definition}\label{def:good_basis}
            A basis $B$ of a lattice $\Lambda$ is called a \textit{good basis} if there exists $T\in \mathcal A_3$ such that $T\Lambda$ is a WR twist with a minimal basis $TB$.
        \end{definition}
    Note that the condition $\alpha\beta\gamma=1$ can be relaxed, since if $\alpha\beta\gamma \ne 1$, we can multiply $\alpha,\beta,\gamma$ by $1/(\alpha\beta\gamma)^{1/3}$, and the resulting lattices after twisting are still similar.

	Now let $B = \{x,y,z\}\subseteq\R^3$ be a basis of a lattice $\Lambda$, and write $x = (x_1,x_2,x_3)$ and similarly for $y$ and $z$.   We wish to compute a matrix $T\in\mathcal{A}_3$ such that $T\Lambda$ is WR.  We can force this condition by insisting that the twisted basis $TB = \{Tx, Ty, Tz\}$ is a minimal basis of the \emph{twist} $T\Lambda$, in the sense of \Cref{sec:wr_lattice}. In particular, we must have $\|Tx\|^2 = \|Ty\|^2 = \|Tz\|^2$, which we will see determines $T$ up to a scalar multiple.
	
	For any twisting matrix $T$, we have
	\begin{align}\label{eq:twisted_length}
		\|Tx\|^2 &= \alpha^2x_1^2 + \beta^2x_2^2 + \gamma^2x_3^2, \\
		\|Ty\|^2 &= \alpha^2y_1^2 + \beta^2y_2^2 + \gamma^2y_3^2 ,\\
		\|Tz\|^2 &= \alpha^2z_1^2 + \beta^2z_2^2 + \gamma^2z_3^2,
	\end{align}	and	
    \begin{align}\label{eq:twisted_innerprod}
		\langle Tx, Ty \rangle &= \alpha^2 x_1y_1 + \beta^2 x_2y_2 + \gamma^2 x_3y_3 ,\\
		\langle Tx, Tz \rangle &= \alpha^2 x_1z_1 + \beta^2 x_2z_2 + \gamma^2 x_3z_3 ,\\
		\langle Ty, Tz \rangle &= \alpha^2 y_1z_1 + \beta^2 y_2z_2 + \gamma^2 y_3z_3.
	\end{align}
	Verifying that the twisted basis $TB$ is a minimal basis of WR lattice requires us to understand the Gram matrix $G_{TB}$, whose entries are given by the above expressions.

	To begin explicitly computing the twisting matrix $T$ which equalizes the Euclidean norms of the basis vectors, we write out the equations $\|Tx\|^2 = \|Ty\|^2$ and $\|Ty\|^2 = \|Tz\|^2$ in terms of the coordinates of $x$, $y$, and $z$, and the entries of $T$.  This gives us the following linear system in the variables $\alpha^2$, $\beta^2$, and $\gamma^2$:
	\[
	\begin{bmatrix}
		x_1^2-y_1^2 & x_2^2-y_2^2 & x_3^2-y_3^2 \\
		y_1^2-z_1^2 & y_2^2-z_2^2 & y_3^2-z_3^2
	\end{bmatrix}
	\begin{bmatrix}
		\alpha^2 \\
		\beta^2 \\
		\gamma^2
	\end{bmatrix}
	=
	\begin{bmatrix}
		0 \\ 0
	\end{bmatrix}. 
	\]
	The general solution to this linear system is given by 
	\begin{align}\label{eq:alp_beta_gam}
		\alpha^2=k\alpha_0,\beta^2=k\beta_0,\gamma^2=k\gamma_0\ (k\in \R),
	\end{align}
	whereas,	\begin{align}\label{eq:initial_alphabeta}
		\alpha_0 =\begin{vmatrix} 
			x_2^2 - y_2^2 & x_3^2 - y_3^2 \\
			y_2^2 - z_2^2 & y_3^2 - z_3^2
		\end{vmatrix} ,
		\beta_0 = \begin{vmatrix}
			x_3^2 - y_3^2 & x_1^2 - y_1^2 \\
			y_3^2 - z_3^2 & y_1^2 - z_1^2
		\end{vmatrix} ,
		\gamma_0=\begin{vmatrix}
			x_1^2 - y_1^2 & x_2^2 - y_2^2 \\
			y_1^2 - z_1^2 & y_2^2 - z_2^2
		\end{vmatrix}. 
	\end{align} 
To have a solution $k\alpha_0>0,k\beta_0>0,k\gamma_0>0$,  the three determinants in \eqref{eq:initial_alphabeta} must have the same sign. Then one obtains the following proposition.

\begin{proposition}
\label{rem:alphapositive}
Let $\alpha_0,\beta_0,\gamma_0$ as defined in \eqref{eq:initial_alphabeta}. Assume that they have the same sign, and let 
$T_B=\diag{\sqrt{|\alpha_0|},\sqrt{|\beta_0|},\sqrt{|\gamma_0|}}$. Then $B$ is a good basis if and only if the Gram matrix of the basis $T_B B$ is of the form as in \Cref{thm:wr3_nec_suff}.
%and $\mathcal{M}_B$ is the Gram matrix of $T_BB$.  thus obtain $\alpha$, $\beta$, and $\gamma$ as the positive square roots of $k \alpha_0, k \beta_0$, and $k \gamma_0$ respectively.  This puts a condition on our basis $B$ and completely defines our twisting matrix $T$. 	
    
\end{proposition}

The following lemma provides a concise and effective criterion for determining whether any three real numbers have the same sign. This is particularly helpful in our setting, since it preserves the symmetry among $x, y, z$, whereas the individual expressions for $\alpha_0, \beta_0, \gamma_0$ in~\eqref{eq:initial_alphabeta} are not symmetric. To avoid conflict with the notation $x, y, z$ used for vectors, we use $r, s, t$ in the following symmetric criterion.

% \begin{lemma}
% \label{lem:sign_of_three_numbers}
% Let $r, s, t$ be three real numbers. Then
% \begin{enumerate}
%     \item[i)] These numbers are positive if and only if $r + s + t > 0$, $rs + st + tr > 0$, and $rst > 0$.
%     \item[ii)] These numbers are negative if and only if $r + s + t < 0$, $rs + st + tr > 0$, and $rst < 0$.
% \end{enumerate}
% \end{lemma}

% \begin{proof}
% This is obvious for the ``if'' statements. When $r + s + t > 0$, $rs + st + tr > 0$, and $rst > 0$, it implies that either all three are positive, or exactly one of them is positive. If it is the latter case, without loss of generality, we may assume $r > 0$ and $s, t < 0$. This leads to $r > -(s + t)$ and $st > -r(s + t)$, which implies $st > (s + t)^2$, a contradiction. Thus, all three numbers must be positive. A similar argument gives (ii).\\
% \end{proof}
\begin{lemma}
\label{lem:sign_of_three_numbers}
    Let $r,s,t$ be three real numbers and 
        $e_1 = r+s+t,
        e_2= rs + st + tr,
        e_3 =rst.$
    Then $r,s,t$ have the same sign if and only if $e_1e_3$ and $e_2$ are both positive.
\end{lemma}
\begin{proof}
    If $e_1e_3$ is positive, then both $e_1$ and $e_3$ have the same sign. When $r + s + t > 0$, $rs + st + tr > 0$, and $rst > 0$, it implies that either all three are positive, or exactly one of them is positive. If it is the latter case, without loss of generality, we may assume $r > 0$ and $s, t < 0$. This leads to $r > -(s + t)$ and $st > -r(s + t)$, which implies $st > (s + t)^2$, a contradiction. Thus, all three numbers must be positive. A similar argument can use for remaining case.\\
\end{proof}

	% Given the above values of $\alpha^2$, $\beta^2$, and $\gamma^2$, it is a (lengthy) algebra exercise to compute that
	% \begin{align}\label{eq:twisted_length}
	% 	\|Tx\|^2 = \|Ty\|^2 = \|Tz\|^2 = 
	% 	\begin{vmatrix}
	% 		x_1^2 & y_1^2 & z_1^2 \\
	% 		x_2^2 & y_2^2 & z_2^2 \\
	% 		x_3^2 & y_3^2 & z_3^2
	% 	\end{vmatrix}.
	% \end{align}
	
	% The inner product $\langle Tx,Ty\rangle$ is given by
	% \begin{align}\label{eq:twisted_innerprod}
	% 	\langle Tx, Ty \rangle &= \alpha^2 x_1y_1 + \beta^2 x_2y_2 + \gamma^2 x_3y_3 = \begin{vmatrix}
	% 		x_1y_1 & x_2y_2 & x_3y_3 \\
	% 		x_1^2-y_1^2 & x_2^2-y_2^2 & x_3^2-y_3^2 \\
	% 		y_1^2-z_1^2 & y_2^2-z_2^2 & y_3^2-z_3^2
	% 	\end{vmatrix} \\
	% 	&= \begin{vmatrix}
	% 		x_1y_1 & y_1^2 & z_1^2 \\
	% 		x_2y_2 & y_2^2 & z_2^2 \\
	% 		x_3y_3 & y_3^2 & z_3^2 
	% 	\end{vmatrix}
	% 	+
	% 	\begin{vmatrix}
	% 		x_1^2 & x_1y_1 & z_1^2 \\
	% 		x_2^2 & x_2y_2 & z_2^2 \\
	% 		x_3^2 & x_3y_3 & z_3^2
	% 	\end{vmatrix}
	% 	+
	% 	\begin{vmatrix}
	% 		x_1^2 & y_1^2 & x_1y_1 \\
	% 		x_2^2 & y_2^2 & x_2y_2 \\
	% 		x_3^2 & y_3^2 & x_3y_3
	% 	\end{vmatrix}
	% \end{align}
	% and similarly for $\langle Tx,Tz\rangle$ and $\langle Ty,Tz\rangle$.

\subsection{Cyclic cubic fields and its  ideal lattices }\label{sec:cyclic_cubic_fields}

%First, we recall the definition of an ideal lattice in a general number field.

  %  Let $F$ be a number field of degree $n$ and signature $(r_1, r_2)$, where $r_1$ is the number of real embeddings and $2r_2$ is the number of complex embeddings. Then $F$ has exactly $n = r_1 + 2r_2$ distinct embeddings into $\C$, and we may select a set of embeddings up to complex conjugation, denoted $\sigma_1, \dots,  \sigma_{r_1+r_2}$ where the first $r_1$ of them are real, and the last $r_2$ of them are complex. We denote the canonical embedding by $\Phi: F \hookrightarrow F \otimes \mathbb{R} \cong \mathbb{R}^{r_1} \times \mathbb{C}^{r_2}$ the map defined by $\Phi(f)= (\sigma_1(f), \cdots, \sigma_{r_1+r_2}(f))$. Here $\mathbb{R}^{r_1} \times \mathbb{C}^{r_2}$ is a Euclidean space equipped with the scalar product: $\langle u, v\rangle = \sum_{i=1}^{r_1} u_i v_i + 2\sum_{i= r_1+1}^{r_1+r_2} \Re(u_i \overline{v_i}) $ where $\overline{v_i}$ is the complex conjugate of $v_i$. \\
%Let  $Q$ be a (fractional) ideal of $F$. Then it is known that $\Phi(Q)$ is a full-rank lattice in $\mathbb{R}^{r_1} \times \mathbb{C}^{r_2}$  by \cite{Bayer-Fluckiger99}. By identifying  $Q$ and $\Phi(Q)$, one has that  $Q$ is an ideal of $F$ and also a lattice in $\mathbb{R}^{r_1} \times \mathbb{C}^{r_2}$. Hence we call ideals of $F$ ideal lattices, see  \cite{Bayer-Fluckiger99} and also \cite[Section 4]{Schoof08} for more details. An ideal lattice $Q$ is called WR if the lattice  $\Phi(Q)$ is WR. 

   In this section, we focus on the case where $F$ is a cyclic cubic field with conductor $m$. We develop a method to test whether a given basis of an arbitrary (fractional) ideal of $F$ is good or not (see \Cref{rem:twist_ideal_cyclic}). Furthermore, we set up some conditions for certain ideals of $F$ to be principal (see Lemmata \ref{lem:dif_idealPID}	and \ref{lem:I_square_prin}). Finally, we prove that the existence of one good basis leads to infinitely many other good bases by multiplying unit elements (see \Cref{prop:inf_good_bas}).
   
   By \cite[pp.6-10]{maki2006determination}, one has \begin{align}
		\label{conductor}m= \frac{a^2+3b^2}{4}
	\end{align} where $a$ and $b$ are integers satisfying one of
	the following conditions,
	\begin{align} \label{eq:a_b_cubic}
		&a \equiv 2 \mod 3, b \equiv 0 \mod 3\text{ and } b>0 \text{  for  } 3 \not| m, \text{ and}\\\nonumber
		&a \equiv 6 \mod 9, b \equiv 3 \text{  or  } 6 \mod 9 \text{ and } b>0 \text{  for  } 3 | m.
	\end{align}
	
	 Moreover, from \cite{maki2006determination}, the following polynomial,  denoted by $df$, can be used to define $F$,
	\begin{align}\label{df-polynomial-cubic}
		df(x)=\left\{\begin{matrix}x^3 - x^2 + \frac{1-m}{3}x -\frac{m(a-3)+1}{27},& \text{if} \ 3 \not | \ m, \\x^3 -\frac{m}{3}x -\frac{a m}{27},& \text{if}\  3 |m.
		\end{matrix}\right.
	\end{align}
When $3\nmid m$, $\mathcal{O}_F$ has a basis $\{\rho,\sigma(\rho),\sigma^2(\rho)\}$  and when $3\mid m$, $\mathcal{O}_F$ has a basis $\{1,\rho,\sigma(\rho)\}$ with $\rho$ as a root of $df(x)$.

 Here we recall the definition of an ideal lattice in $F$. Let $\sigma$ be a generator of the Galois group $G= \Gal(F/\Q)$. Then $F$ has three distinct embeddings into $\R$: $1, \sigma, \sigma^2$.  We denote the canonical embedding by $\Phi: F \hookrightarrow F \otimes \R \cong  \R^{3} $ where $\Phi(f)= (f, \sigma(f), \sigma^2(f))$ for $f\in F$. Here $\R^3$ is a Euclidean space equipped with the usual scalar product: $\langle u, v\rangle = \sum_{i=1}^3 u_i v_i $ for $u=(u_1, u_2, u_3), v= (v_1, v_2, v_3) \in \R^3$.

Let  $\mathfrak Q$ be a (fractional) ideal of $F$. Then it is known that $\Phi(\mathfrak Q)$ is a full-rank lattice in $\R^3$  by \cite{Bayer-Fluckiger99}. By identifying  $\mathfrak Q$ and $\Phi(\mathfrak Q)$, one has that $\mathfrak Q$ is an ideal of $F$ and also a lattice in $\mathbb{R}^3$. Hence, we call ideals of $F$ ideal lattices, see  \cite{Bayer-Fluckiger99} and also \cite[Section 4]{Schoof08} for more details. 
\begin{definition}
    A (fractional) ideal  $\mathfrak Q$ of $F$ is called a \textit{WR ideal (lattice)} if the lattice  $\Phi(\mathfrak Q)$ in $\mathbb{R}^3$ is WR. 
\end{definition}
%%%%%%%%%%%%Add
Let $\mathfrak I$ be an ideal lattice of $\mathcal{O}_F$ and $B=\set{x,y,z}$ be an integral basis of $\mathfrak I$. Note that $x=\circbrack{x,\sigma(x),\sigma^2(x)}$ and similarly for $y$ and $z$. 
 In this case, the formulae for computing $\alpha_0,\beta_0,\gamma_0$ in \eqref{eq:initial_alphabeta} can be rewritten as below. 
\begin{align}\label{eq:alpha2_cylic}
		\alpha_0= \begin{vmatrix}
			\sigma(x)^2 - \sigma(y)^2 & \sigma^2(x)^2 - \sigma^2(y)^2 \\
			\sigma(y)^2 - \sigma(z)^2 & \sigma^2(y)^2 - \sigma^2(z)^2
		\end{vmatrix} ,
		\beta_0 = \begin{vmatrix}
			\sigma^2(x)^2 - \sigma^2(y)^2 & x^2 - y^2 \\
			\sigma^2(y)^2 - \sigma^2(z)^2 & y^2 - z^2
		\end{vmatrix} ,
		\gamma_0 = \begin{vmatrix}
			x^2 - y^2 & \sigma(x)^2 - \sigma(y)^2 \\
			y^2 - z^2 & \sigma(y)^2 - \sigma(z)^2
		\end{vmatrix}
	\end{align}  
One has
\begin{align}    \label{eq:al_0_cycliccu}\beta_0=\sigma(\alpha_0), \gamma_0=\sigma^2(\alpha_0)
\end{align} and thus, three elementary symmetric polynomials in $\alpha_0,\beta_0,\gamma_0$ belong to $\Z$. In other words, \Cref{lem:sign_of_three_numbers} is an effective tool to check whether $\alpha_0,\beta_0,\gamma_0$ have the same sign or not.

\begin{algorithm}[H]
\caption{Testing good bases}\label{rem:twist_ideal_cyclic}
\KwIn{A $\mathbb Z$-basis $B = \{x,y,z\}$ of a fractional ideal $\mathfrak J$ of a cyclic cubic field $F$}
\KwOut{$B$ is good or not}

 \begin{enumerate}
        \item Compute  $\alpha_0,\beta_0,\gamma_0$ as in \eqref{eq:alpha2_cylic}
                and denote by \begin{align*}           e_{1B}&=\alpha_0+\beta_0+\gamma_0,\\
            e_{2B}&=\alpha_0\beta_0+\beta_0\gamma_0+\gamma_0\alpha_0,\\
            e_{3B}&=\alpha_0\beta_0\gamma_0.
        \end{align*}
        \item If $e_{1B}\cdot e_{3B}\le 0$ or $e_{2B}\le 0$, then $B$ is not good by \Cref{lem:sign_of_three_numbers}. Otherwise, move to step 3. 
        \item Compute the matrix $\mathcal{M}_B=\Matrix{s_B&u_B&v_B\\u_B&s_B&w_B\\v_B&w_B&s_B}$ as in \eqref{eq:twisted_length} and \eqref{eq:twisted_innerprod}, \begin{align*}
            s_B&= |\alpha_0|x^2 +|\beta_0|\sigma(x^2)+|\gamma_0|\sigma^2(x^2),\\
            u_B&=|\alpha_0|xy+|\beta_0|\sigma(xy)+|\gamma_0|\sigma^2(xy),\\
            v_B&=|\alpha_0|xz+|\beta_0|\sigma(xz)+|\gamma_0|\sigma^2(xz),\\
            w_B&=|\alpha_0|yz+|\beta_0|\sigma(yz)+|\gamma_0|\sigma^2(yz).
        \end{align*}
        \item Check whether the matrix  $\mathcal{M}_B$ is of the form in  \Cref{thm:wr3_nec_suff} or not. If it is, then  $B$ is a good basis, $T_B=\diag{\sqrt{|\alpha_0|},\sqrt{|\beta_0|},\sqrt{|\gamma_0|}}$ and $\mathcal{M}_B$ is the Gram matrix of the basis $T_BB$ (of the twist lattice).
    \end{enumerate}
\end{algorithm}
\begin{remark}

\begin{enumerate}
\item  \Cref{rem:twist_ideal_cyclic} works for all full-rank lattices in $\mathbb{R}^3$ using $\alpha_0, \beta_0, \gamma_0$ in \ref{eq:initial_alphabeta}. In case this lattice is an ideal of a cyclic cubic field, the entries of $\mathcal{M}_B$ in step 3 of \Cref{rem:twist_ideal_cyclic}  are rational numbers since $$s_B=\Tr\circbrack{|\alpha_0| x^2}, u_B=  \Tr\circbrack{|\alpha_0|xy}, v_B=\Tr\circbrack{|\alpha_0|xz},w_B=\Tr\circbrack{|\alpha_0|yz}.$$

\item In case the ideal is the ring of integers $\mathcal{O}_F$, and assume that it has a WR twist with a good basis $B=\{x,y,z\}$, then all $\alpha_0,\beta_0,\gamma_0$ (obtained by applying \eqref{eq:alpha2_cylic}) have the same sign and the matrix $\mathcal{M}_B$ must be of the form as in \Cref{thm:wr3_nec_suff}.

\end{enumerate}

\end{remark}

Note that the discriminant of $F$ is $\Delta_F=m^2$.
	Then, by \cite[Theorem 4.8]{conrad2009different}, the different ideal of $F$ has norm $m^2$.
    
	\begin{lemma}\label{lem:dif_idealPID}
		Let $F$ be a cyclic cubic field of the conductor $m$ defined as in the beginning of \Cref{sec:cyclic_cubic_fields}. Then the different ideal of $F$ is principal  if $b=3$, or $m$ is prime.
      	\end{lemma}
	\begin{proof}
		In the first case, when $b=3$, we can verify that the absolute value of the norm of \( df'(\rho) \) is \( m^2 \), here $df(x)$ is the defining polynomial of $F$ as in \eqref{df-polynomial-cubic} and $df'$ is the derivative of $df$. Indeed,
		\begin{itemize}
			\item If \( 3 \nmid m \), then \( df'(\rho) = 3\rho^2 - 2\rho + \dfrac{1 - m}{3} \), and a direct computation in \texttt{Mathematica} shows that  \( \Nm(df'(\rho)) = -m^2 \).
			\item If \( 3 \mid m \), then \( df'(\rho) = 3\rho^2 - \dfrac{m}{3} \), and then \( \Nm(df'(\rho)) = m^2 \).
		\end{itemize}
		Since the different ideal has norm $m^2$, and the ideal of norm \( m^2 \) is unique, it follows that the different ideal of \( F \) is generated by \(df'(\rho)\).

        Now consider the second case,   $m$ is prime. Then the Sylow $3-$subgroup of the class group of $F$ must be trivial as the result of \cite{gras1973ell} and \cite[Theorem 1]{mayer2022theoretical}. Since the norm of the different ideal is $m^2$, its equivalent class in the class group $\text{Cl}(F)$ belongs to the Sylow $3-$subgroup. It means that this class is the identity of $\text{Cl}(F)$. Therefore, the different ideal is principal.
	\end{proof}
	
	\begin{lemma}\label{lem:squareideal}
		Let $F$ be a cyclic cubic field and $\mathfrak I$ be an ideal of $\mathcal{O}_F$ such that $\Nm(\mathfrak I)\mid \Delta_F^t$ for some $t\in \Z_{>0}$. Then there exists $q\in \Q_{>0} $ and such that $\mathfrak I=q\mathfrak I^4$.
	\end{lemma}
	\begin{proof}
		We have $\mathfrak I^3=\Nm(\mathfrak I)\mathcal{O}_F$. It implies $\Nm(\mathfrak I)\mathfrak I=\mathfrak I^4 $. Therefore, $\mathfrak I = \dfrac{1}{\Nm(\mathfrak I)}\mathfrak I^4$. 
	\end{proof}

\begin{lemma}\label{lem:I_square_prin}
    Let $F$ be a cyclic cubic field and $\mathfrak I$ be an ideal of $\mathcal{O}_F$ that $N(\mathfrak I)\mid \Delta_F^t$ for some $t\in \Z_{>0}$. Then $\mathfrak I$ is principle if and only if $\mathfrak I^2$ is principle.
\end{lemma}
\begin{proof}
Let $n=N(\mathfrak I).$ Since $n\mid\Delta_F^t$, then $p\mathcal{O}_F=\mathfrak P^3$ for all primes $p\mid n$ and all unique prime ideal $\mathfrak P$ of $\mathcal{O}_F$ over $p$. It is easy to  verify that $n\mathcal{O}_F= \mathfrak I^3.$ If $\mathfrak I$ is principle then $\mathfrak I^2$ is obviously principle. If $\mathfrak I^2$ is principle then $\mathfrak I^4$ is principle and we write $\mathfrak I^4=(a)$ for some $a\in \mathfrak I^4$. It implies that $\mathfrak I$ is principle ideal as $\mathfrak I^3$ is also principle.
\end{proof}
Once $\mathcal{O}_F$ admits a WR twist $L$, there are infinitely many bases of $\mathcal{O}_F$ that yield twists similar to $L$. 
\begin{proposition}\label{prop:inf_good_bas}
    Let $\mathfrak{I}$ be an ideal of $\mathcal{O}_F$ and $u$ be a unit. If $\mathfrak{I}$ has a good basis $B=\set{x,y,z}$, then $B'=\set{ux,uy,uz}$ is also a good basis of $\mathfrak{I}$. Moreover, $B$ and $B'$ yield the same twist lattice (up to similarity).
\end{proposition}
\begin{proof}
    It is clear that $B'$ is also a basis of $\mathfrak{I}$. We compute the quantities $\alpha_0$, $\beta_0$, and $\gamma_0$ for $B$ as in \eqref{eq:alpha2_cylic}.

    The Gram matrix $G_B$ of $T_BB$ is \begin{align*}
        G_B=\Matrix{s_B&u_B&v_B\\u_B&s_B&w_B\\v_B&w_B&s_B}
    \end{align*}where \begin{align*}
    s_B& = \alpha_0 x^2+\beta_0\sigma(x^2)+\gamma_0 \sigma^2(x^2),\\
        u_B&= \alpha_0 xy+\beta_0\sigma(xy)+\gamma_0 \sigma^2(xy),\\
        v_B&= \alpha_0 xz+\beta_0\sigma(xz)+\gamma_0 \sigma^2(xz),\\
        w_B&= \alpha_0 yz+\beta_0\sigma(yz)+\gamma_0 \sigma^2(yz).
    \end{align*}
    By \Cref{rem:twist_ideal_cyclic}, $\alpha_0,\beta_0,\gamma_0$ have the same sign. Let $\alpha_0'=\dfrac{\alpha_0}{u^2},\beta_0'=\dfrac{\beta_0}{\sigma(u^2)},\gamma_0'=\dfrac{\gamma_0}{\sigma^2(u^2)}$. Then $\alpha_0',\beta_0',\gamma_0'$ have the same sign. Let $T'_B=\diag{\sqrt{|\alpha_0'|},\sqrt{|\beta_0'|},\sqrt{|\gamma_0'|}}$. Then the Gram matrix $G_{B'}$ of $T_{B'}B'$ is \begin{align*}
    G_{B'}=\Matrix{s'_{1B}&u'_B&v'_B\\u'_B&s'_{2B}&w'_B\\v'_B&w'_B&s'_{3B}}
    \end{align*}
    It is straightforward to prove that $G_B=G_{B'}$ by calculation.

\end{proof}
    \section{Well-roundedness of ideals of norm dividing the discrimiant}\label{sec:ideal_norm_div_dis}
From now on, we denote by $F$  a cyclic cubic field and $\mathcal{O}_F$ its ring of integers. Let $\sigma$, $m$ and $\Delta_F$ be a generator of $\Gal(F/\Q)$, the  conductor and discriminant of $F$, respectively.

    In this section,  we first compute all integral bases of ideal lattices $\mathfrak I$ whose norms divide $\Delta_F$.  We then prove that if such an ideal $\mathfrak I$ is WR, then it must have a minimal basis of the form $\set{\pm \kappa,\pm\sigma(\kappa),\pm\sigma^2(\kappa)}$  where $\kappa$ is a shortest vector of $\mathfrak I$.  Note that the result of this section will be applied to prove the form of the Gram matrix for a basis that satisfies the conditions of \Cref{thm:prin_ideal_from_WRtwist}. 

    Many WR ideals with norm dividing the discriminant in a cyclic cubic field are described in  ~\cite{tran2023well}, it suffices to complete the remaining cases on WR ideals that were not settled in this work. In other words, with the notation defined below and by \Cref{rem:done_ideals}, we only consider ideals of the form $\mathfrak P_I^2\mathfrak P_J$ for the case $3\nmid m$ and ideals of the form $\mathfrak P_I^2\mathfrak P_J,\mathfrak P_0^2\mathfrak P_I,\mathfrak P_0^2\mathfrak P_I^2\mathfrak P_J$ for the case $3\mid m$.

The conductor $m$ can be written as $m=p_1p_2\cdots p_r$ if $3\nmid m$ and $m=9p_1p_2\cdots p_r$ if $3\mid m$ where $p_i$ are prime and $p_i\equiv 1\pmod 3$. Let $I,J$ be subsets of the set $\set{1,2,\cdots, r}$ and $p_I= \prod_{i\in I}p_i, p_J=\prod_{i\in J}p_i$ and $\mathfrak P_I=\prod_{i\in I}\mathfrak P_i,\mathfrak P_J=\prod_{j\in J}\mathfrak P_j$ where $\mathfrak P_i$ is the unique prime ideal above $p_i$.  We also denote $\mathfrak P_0$ as the unique prime ideal above $3$ in the case $3\mid m.$

\begin{remark}\label{rem:done_ideals}
    Some conditions for certain ideals of norm dividing $\Delta_F$ to be WR are presented in \cite{tran2023well}. In particular, for the case $3\nmid m$, the ideal $\mathfrak P_I$ is WR if and only if $\dfrac{m}{4}\le p_I^2\le 4m$. In the other case, the ideal $\mathfrak P_I$ is not WR, the ideal $\mathfrak P_0\mathfrak P_I$ is WR iff and only if $\dfrac{m}{36}\le p_I^2\le \dfrac{4m}{9}$, and the ideal $\mathfrak P_0\mathfrak P_I^2\mathfrak P_J$ is WR if and only if $\dfrac{m}{36}\le p_I^2p_J\le \dfrac{4m}{9}$. When these ideals are WR, they must have a minimal basis of the form $\set{\pm \kappa,\pm\sigma(\kappa),\pm\sigma^2(\kappa)}$. 
    
    We note that the condition $N(\mathfrak I)\mid\Delta_F$ is one of the sufficient conditions for an ideal $\mathfrak I$ to be WR by \cite[Section 4.2] {tran2023well}. An unproven conjecture given in \cite[Section 5] {tran2023well} says that this is also a necessary condition. 
\end{remark}

When $3\nmid m$, let $\beta=\alpha-\dfrac{1}{3}$. Then $\beta,\sigma(\beta), \sigma^2\circbrack{\beta}$ are all roots of the polynomial \begin{align}
    \label{eq:beta}x^3-\dfrac{m}{3}x-\dfrac{am}{27}=0.
\end{align} 
Note that every element $\delta $ of $F$ is of the form $\delta=m_1+m_2\beta+m_3\sigma(\beta)$ for some rational number $m_1,m_2,m_3$ and thus, by applying a similar argument in the proof of \cite[Lemma 3.12]{tran2023well}, one obtains  \begin{align}
    \label{eq:length_fom} \|\delta\|^2=3m_1^2+\dfrac{2}{3}m\circbrack{m_2^2-m_2m_3+m_3^2}.
\end{align}
Without loss of generality, we can choose $\sigma$ such that \begin{align}
    \label{eq:t}\mathfrak{t}=\circbrack{\beta-\sigma(\beta)}\circbrack{\sigma(\beta)-\sigma^2(\beta)}\circbrack{\sigma^2(\beta)-\beta}\le 0.
\end{align}

It leads to the fact that \begin{align}
\nonumber&\psi_1=\beta+\sigma\circbrack{\beta}+\sigma^2\circbrack{\beta}=  0,\\\label{eq:psi123}&\psi_2=\beta\sigma(\beta)+\sigma(\beta)\sigma^2(\beta)+\sigma^2(\beta)=\dfrac{-m}{3},\\\nonumber
&\psi_3=\beta\sigma(\beta)\sigma^2(\beta)=\dfrac{am}{27}.
\end{align} 

\begin{lemma}\label{lem:3nmidm_yIzI}
	There exist integers $y_I,z_I$ such that $y_I+z_I\equiv 1\pmod 3$, $p_I= y_I^2-y_Iz_I+z_I^2$ and $p_I\mid (y_I\circbrack{a+3b}+z_I\circbrack{a-3b)}$.
\end{lemma}
\begin{proof}
	One has $\dfrac{m}{9}=\circbrack{\dfrac{a+b}{6}}^2-\dfrac{a+b}{6}\dfrac{b}{3}+\dfrac{b^2}{9}$. It implies that $m=  \circbrack{\dfrac{a+b}{2}}^2-\dfrac{a+b}{2}b+b^2$. Using a similar argument in the proof of \cite[Lemma 3.23]{tran2023well}, there exists $y_I,z_I$  such that $p_I=y_I^2-y_Iz_I+z_I^2$ and $y_I+z_I\equiv 1\pmod 3$ (we change signs of results) and \begin{align*}
		p_I\mid\circbrack{a+b}y_I+2bz_I-(a+b)z_I,2by_I-(a+b)z_I .
	\end{align*}
    It leads to $p_I\mid \circbrack{a+b}y_I+\circbrack{b-a}z_I,p_I\mid 2by_I-(a+b)z_I $ and we obtain $p_I\mid (z_I\circbrack{3b+a}+y_I\circbrack{a-3b})$.
    
\end{proof}
 
\begin{lemma}
Let $y_I,z_I$ be as in \Cref{lem:3nmidm_yIzI} and $x= \dfrac{p_Ip_J-y_I-z_I}{3}$. Then $x\in \Z$.    If $3\nmid m$ and $I,J$ are disjoint subsets of $\set{1,2,\cdots,r}$, then $\set{\kappa,\sigma(\kappa),\sigma^2(\kappa)}$ is a basis of $\mathfrak P_I^2\mathfrak P_J$ and $\mathfrak P_I^2\mathfrak P_J$ is WR if and only if $\dfrac{m}{4}\le p_I^2p_J\le 4m$ where $\kappa =  x+y_I\alpha +z_I\sigma(\alpha)$ .
\end{lemma}

\begin{proof}
    We first prove that $\kappa \in \mathfrak P_I^2\mathfrak P_J$. Indeed, with the notation as in \eqref{eq:t} and \eqref{eq:psi123}, one has \begin{align*}
    \T^2&= -27\psi_3^2-4\psi_1^3\psi_3+\psi_2^2\psi_1^2+18\psi_3\psi_2\psi_1-4\psi_2^3=\dfrac{4m^3}{27}-\dfrac{a^2m^2}{27}=\circbrack{\dfrac{bm}{3}}^2,\\
        \T&= \circbrack{\beta^2\sigma^2(\beta)+\sigma(\beta)^2\beta+\sigma^2(\beta)^2\sigma(\beta)}-\circbrack{\beta^2\sigma(\beta)+\sigma(\beta)^2\sigma^2(\beta)+\sigma^2(\beta)^2\beta}.
    \end{align*}
    We implies that $\circbrack{\beta^2\sigma^2(\beta)+\sigma(\beta)^2\beta+\sigma^2(\beta)^2\sigma(\beta)}-\circbrack{\beta^2\sigma(\beta)+\sigma(\beta)^2\sigma^2(\beta)+\sigma^2(\beta)^2\beta}=-\dfrac{mb}{3}$. Furthermore, since $\circbrack{\beta^2\sigma^2(\beta)+\sigma(\beta)^2\beta+\sigma^2(\beta)^2\sigma(\beta)}+\circbrack{\beta^2\sigma(\beta)+\sigma(\beta)^2\sigma^2(\beta)+\sigma^2(\beta)^2\beta}+3\psi_3=  \psi_1\psi_2=~0$, we have the following equations\begin{align*}
\beta^2\sigma^2(\beta)+\sigma(\beta)^2\beta+\sigma^2(\beta)^2\sigma(\beta)&= \dfrac{m}{18}\circbrack{-3b-a},\\\beta^2\sigma(\beta)+\sigma(\beta)^2\sigma^2(\beta)+\sigma^2(\beta)^2\beta&=\dfrac{m}{18}\circbrack{-3b+a}.
    \end{align*} 
  Let $X = \dfrac{3x+y_I+z_I}{3}$. Then one can write $\kappa$ as follow \begin{align*}
    \kappa =X+y_I\beta+z_I\sigma(\beta).
\end{align*}

    It is quite straightforward to prove that \begin{align*}
    \Nm\circbrack{\kappa}=& X^3+X^2(y_I+z_I)\psi_1+X\circbrack{y_I^2-y_Iz_I+z_I^2}\psi_2+Xy_Iz_I\psi_1^2\\&+y_1^2z_1\psi_1\psi_2+\circbrack{y_I^3-3y_I^2z_I+z_I^3}\psi_3+\circbrack{y_Iz_I^2-z_Iy_I^2}\circbrack{\beta^2\sigma^2(\beta)+\sigma^2(\beta)\beta+\sigma^2(\beta)\sigma(\beta)}\\ =& X^3-\dfrac{mX}{3}p_I+\circbrack{y_I+z_I}p_I
    \dfrac{am}{27}\\&+y_Iz_I^2\circbrack{\beta^2\sigma^2(\beta)+\sigma(\beta)^2\beta+\sigma^2(\beta)^2\sigma(\beta)}+y_I^2z_I\circbrack{\beta^2\sigma(\beta)+\sigma(\beta)^2\sigma^2(\beta)+\sigma^2(\beta)^2\beta}\\=&\dfrac{p_I^3p_J^3}{27}-\dfrac{mp_Ip_J}{9}+\circbrack{y_I+z_I}p_I\dfrac{am}{27}-\dfrac{m}{18}y_Iz_I\circbrack{z_I\circbrack{3b+a}+y_I\circbrack{a-3b}}.
\end{align*}
Using \Cref{lem:3nmidm_yIzI}, we obtain $ p_I^2p_J\mid 54\Nm\circbrack{\kappa}.$ Combining with $\kappa\in \mathcal{O}_F$, one obtains that the element $\kappa$ belongs to $\mathfrak P_I^2\mathfrak P_J.$ It also leads to $L=\kappa \Z\oplus \sigma(\kappa)\Z\oplus \sigma^2(\kappa)\Z\subset \mathfrak P_I^2\mathfrak P_J$. Since $\textrm{Vol}\circbrack{L}=p_I^2p_J=\Nm\circbrack{\mathfrak P_I^2\mathfrak P_J}$, one implies $L= \mathfrak P_I^2\mathfrak P_J.$

Let 
\begin{align*}
    \delta  &=  m_1 \kappa+m_2\sigma(\kappa)+m_3\sigma^2\circbrack{\kappa}\\&=\circbrack{m_1+m_2+m_3}X+\circbrack{m_1y-m_2z+m_3(z-y)}\beta+\circbrack{m_1z+m_2(y-z)-m_3y}\sigma\circbrack{\beta}.
\end{align*}Then \begin{align*}
    \|\delta\|^2&=3\circbrack{m_1+m_2+m_3}^2X^2+\dfrac{2m}{3}\circbrack{y_I^2-y_Iz_I+z_I^2} \circbrack{m_1^2+m_2^2+m_3^2-m_1m_2-m_1m_3-m_2m_3}.
\end{align*}
By using a similar argument as the one in the proof of \cite[Proposition 3.20]{tran2023well}, one has 
\begin{align*}
	\min_{\delta \ne 0}\|\delta\|^2=\min\set{3p_I^2p_J^2,2mp_I,\dfrac{p_I^2p_J^2}{3}+\dfrac{2m}{3}p_I}
\end{align*}
and the lattice $\mathfrak P_I^2\mathfrak P_J$ is WR if and only if $\min_{\delta \ne 0}\delta =\dfrac{p_I^2p_J^2}{3}+\dfrac{2m}{3}p_I$ . It is equivalent
to the statement
\begin{align*}
	\dfrac{p_I^2p_J^2}{3}+\dfrac{2mp_I}{3}\le 3p_I^2p_J^2,\dfrac{p_I^2p_J^2}{3}+\dfrac{2mp_I}{3}\le 2mp_I.
\end{align*}In other words, $\mathfrak P_I^2\mathfrak P_J$ is WR if and only if $\dfrac{m}{4}\le p_I^2p_J \le 4m$.
    \end{proof}
The set of minimal vectors of the ideal $\mathfrak P_I^2\mathfrak P_J$ is $\set{\pm \kappa,\pm\sigma(\kappa), \pm\sigma^2(\kappa)}$. Combining this fact and  \cite[Proposition 3.8]{tran2023well}, we have the following proposition.
\begin{proposition}\label{prop:con_basis_3nmidm}
	Let $F$ be a cyclic cubic field of which the discriminant $\Delta_F$ is not a multiple of $9$ and let $\mathfrak I$ be an ideal of $\mathcal{O}_F$ such that $\Nm\circbrack{\mathfrak I}\mid \Delta_F$. If $\mathfrak I$ is WR and $\kappa$ is a shortest vector of $\mathfrak I$, then every minimal basis of $\mathfrak I$ is of the form $\set{\pm \kappa, \pm \sigma\circbrack{\kappa}, \pm \sigma^2\circbrack{\kappa}}.$    
\end{proposition}
When $3\mid m$, let $\dfrac{m}{9}=A^2-AB+B^2$ and $x_I$ be a root of the equation $Bx\equiv -A \pmod {p_I}$. 
\begin{lemma}\label{lem:ideal_3midm}
    If $3\mid m$, $p_0=3$ and $I,J$ are disjoint subsets of $\set{1,2,\cdots,r}$, then we have the following.
    \begin{enumerate}
        \item[i)] $\mathfrak P_I^2\mathfrak P_J= p_Ip_J\Z\oplus p_I\alpha \Z\oplus \circbrack{x_I\alpha -\sigma(\alpha)}\Z.$
        \item [ii)]$\mathfrak P_0^2\mathfrak P_I=p_0p_I\Z\oplus3\alpha\Z \oplus\circbrack{\alpha-\sigma(\alpha)}\Z$.
        \item[iii)]$\mathfrak P_0^2\mathfrak P_I^2\mathfrak P_J= p_0p_Ip_J\Z \oplus p_0p_I\alpha \Z \oplus \circbrack{x\alpha-\sigma(\alpha)}\Z$.
    \end{enumerate}
\end{lemma}
\begin{proof}
By calculating the norms, it is easy to prove that the lattices on the right-hand side are contained in those on the left-hand side. Moreover, since every pair of these left and right-hand side lattices has the same covolume, they must be identical. 
\end{proof}
\begin{lemma}\label{lem:lem_not_WR}
    Let $3\mid m$ and $L$ be a full-rank sublattice of $\mathcal{O}_F$. If $L$ has an integral basis of the form $\set{k,\kappa_L=x_L\alpha+y_L\sigma(\alpha), \tau_L=u_L\alpha+v_L\sigma\circbrack{\alpha}}$ where $k, x_L,y_L,u_L,v_L\in \Z$, then $L$ is not WR.
\end{lemma}
    \begin{proof}
        Let $M_2= m_2x_L+m_3u_L$, $M_3=m_2y_L+m_3v_L$ and  $\delta = m_1k+m_2\kappa_L+m_3\tau_L$ where $m_1,m_2,m_3\in \Z$. Then by \cite[Lemma 3.12]{tran2023well}, we have \begin{align*}
            \|\delta\|^2 =  3m_1^2k^2+\dfrac{2m}{3}\circbrack{M_2^2-M_2M_3+M_3^2}.
        \end{align*}
        If $m_1\ne 0$ and $ M_2^2-M_2M_3+M_3^2\ne 0$, then $\|\delta\|^2\ge 3m_1^2k^2\ge 3k^2= \| k\|^2$. Note that $M_2^2-M_2M_3+M_3^2=0$ if and only if $M_2=M_3=0$. We imply that \begin{align*}
            \min_{\delta\in L^*}\|\delta\|^2=\min\set{3k^2, \min_{\delta\in L\cap \langle \alpha,\sigma(\alpha)\rangle}\|\delta\|^2}.
        \end{align*} The numbers $3k^2$ and $\min_{\delta\in L\cap \langle \alpha,\sigma(\alpha)\rangle}\|\delta\|^2$ are different. Indeed, suppose, for the sake of contradiction, that $3k^2=\min_{\delta\in L\cap \langle \alpha,\sigma(\alpha)\rangle}\|\delta\|^2$. There exists $g,h\in \Z$ such that \begin{align}\label{eq:eq1}k^2=2p_1\cdots p_r\circbrack{g^2-gh+h^2}.
        \end{align} Let $v_3\circbrack{k^2}=\max\set{i\in \mathbb{N}:3^i\mid k^2}$. Then $v_3\circbrack{k^2}$ is even. Dividing both sides of \eqref{eq:eq1} by \(3^{v_3(k^2)}\), we obtain that the left-hand side is congruent to $1$ modulo $3$, while the right-hand side is congruent to $2$ modulo $3$, which is a contradiction. Therefore, the minimal vectors belong to one of the two sets $ L\cap \Z$ or $L\cap \langle \alpha, \sigma(\alpha)\rangle$ but not both (here $\langle \alpha, \sigma(\alpha)\rangle= \{ m_1 \alpha + m_2 \sigma(\alpha): m_1, m_2 \in \mathbb{Z}\}$). However, in these cases, there are at most $2$ independent shortest vectors in each set. Consequently, $L$ is not WR.
    \end{proof}
    The following lemma is a consequence of \Cref{lem:ideal_3midm} and \Cref{lem:lem_not_WR}.
\begin{lemma}\label{lem:ideal_3midm_not_WR}
      If $3\mid m$, $p_0=3$ and $I,J$ are disjoint subsets of $\set{1,2,\cdots,r}$. Then the ideals $\mathfrak P_I^2\mathfrak P_J,\mathfrak P_0^2\mathfrak P_I,\mathfrak P_0^2\mathfrak P_I^2\mathfrak P_J$ are not WR.
\end{lemma}

Combining \Cref{lem:ideal_3midm_not_WR} with \cite[Propositions 3.20 and 3.25]{tran2023well}, we obtain the following proposition.
\begin{proposition}\label{prop:con_basis_3midm}
    Let $F$ be a cyclic cubic field of which the discriminant $\Delta_F$ is a multiple of $9$, and let $\mathfrak I$ be an ideal of $\mathcal{O}_F$ such that $\Nm\circbrack{\mathfrak I}\mid \Delta_F$. If $\mathfrak I$ is WR and $\kappa$ is a shortest vector of $\mathfrak I$, then every minimal basis of $\mathfrak I$ is of the form $\set{\pm \kappa, \pm \sigma\circbrack{\kappa}, \pm \sigma^2\circbrack{\kappa}}.$\end{proposition}
% As a consequence of \Cref{prop:con_basis_3midm} and \Cref{prop:con_basis_3nmidm}, the Gram matrix of every minimal basis of the ideal $I$ in these propositions is of the form as below. 
    \begin{corollary}\label{cor:Gram_matrix}
        Let $\mathcal{O}_F$ be the ring of integers of a cyclic cubic field with discriminant $\Delta_F$ and $\mathfrak I$ be a WR ideal of $\mathcal{O}_F$ such that $\Nm\circbrack{\mathfrak I}\mid \Delta^t_F$ for some positive integer $t$. Then the Gram matrix of every minimal basis of $\mathfrak I$ is of the form  \begin{align*}
        \begin{bmatrix}
			s & u & v \\
			u & s & w \\
			v & w & s
		\end{bmatrix}
    \end{align*}
    where $|u|=|v|=|w|.$

    \begin{proof}
        The statement is true when $t=1$ as a consequence of \Cref{prop:con_basis_3midm,prop:con_basis_3nmidm}. In the case $t>1$, since $\Nm(\mathfrak{I})\mid \Delta_F^t$, all prime factors of $\Nm(\mathfrak{I})$ are ramified, thus for each prime $p\mid \Nm(\mathfrak{I})$, $\mathfrak{P}^3= p\mathcal{O}_F$ where $\mathfrak{P}$ is a unique prime ideal above $p$. One implies that there exists an ideal $\mathfrak{I}$ and an integer $q>0$ such that $\mathfrak{I}=q\mathfrak{J}$ and $\Nm\circbrack{\mathfrak{J}}\mid \Delta _F$. In other words, $\mathfrak{I}$ is similar to $\mathfrak{J}$. The result is implied by applying the case $t=1$.
    \end{proof}
        
    \end{corollary}
 
	%%%%%%%%%%%%%%%%%%%%%%%%%%%%%%
	\section{WR twists for the ring of integers of cyclic cubic fields}\label{sec:wr_twist_of_roi}
	In this section, we investigate conditions for the ring of integers of a cyclic cubic field $F$ to have a WR twist.    
    It is known that for any lattice with an arbitrary rank, it always has a WR twist as a result of Solan \cite{solan2019stable}. In other words, a WR twist of $\mathcal{O}_F$ always exists. However, Solan's proof does not show us how to compute this twist explicitly.  
    
    We will show that if the ring of integers $\mathcal{O}_F$ has a WR and principal ideal, then one can construct a WR twist for $\mathcal{O}_F$ that is similar to this ideal lattice. Conversely, if $\mathcal{O}_F$ has a WR twist satisfying additional conditions (see \Cref{thm:prin_ideal_from_WRtwist}), then there exists a principal ideal of $\mathcal{O}_F$ that is similar to the WR twist from $\mathcal{O}_F$ and the Gram matrix of the corresponding twisted basis must have equal off-diagonal entries in absolute value.

     Additionally, we will show later in this section that $\mathcal{O}_F$ has an orthogonal WR twist when the different ideal of $F$ is principal and generated by a totally positive element (see \Cref{thm:ortho_twist} and \Cref{prop:cyc_cu_orth}).

\begin{proposition}\label{prop:WRtwist_from_prin_ideal}
Let $F$ be a cyclic cubic field. Suppose that the ring of integers $\mathcal{O}_F$ has a WR and principal ideal $\mathfrak I$ with a minimal basis $ \{a, b,c \}$. Then $\mathfrak I$ is a WR twist of  $\mathcal{O}_F$  up to a similarity.
\end{proposition}
\begin{proof} Let $\delta$ be a generator of $\mathfrak I$, i.e., $\mathfrak I=\langle \delta\rangle$.
		
Let $B=\set{x,y,z}$ be a $\Z$-basis for $\mathcal{O}_F$. Then $B_\mathfrak I=\set{\delta x,\delta y,\delta z}$ is a $\Z$-basis for $\mathfrak I$. Let $B_{\mathfrak I'}=\{a,b,c\}$ be a minimal basis for $\mathfrak I$. Then there exists a unimodular matrix $U$ such that 
\begin{align*}
			\left[ \begin{matrix}
				\delta x&\delta y&\delta z
			\end{matrix}\right]U = \left[\begin{matrix}
				a & b & c
			\end{matrix}\right].
		\end{align*}
		One has 
		\begin{align*}
			G_{B_{\mathfrak I'}}=B_{\mathfrak I'}^TB_{\mathfrak I'}=U^TB^T\diag{\delta^2,\sigma(\delta^2),\sigma^2(\delta^2)}BU
		\end{align*}
		where $\sigma$ is a generator of the Galois group of $F/\Q$. Let $B'=BU=\set{x',y',z'}$ and 
		\begin{align*}
			T = \begin{bmatrix}
				\abs{\delta}&0&0\\
				0&\abs{\sigma\circbrack{\delta}}&0\\0&0&\abs{\sigma^2\circbrack{\delta}}
			\end{bmatrix}.
		\end{align*}
		Since $U$ is a  unimodular matrix, $\set{x',y',z'}$ is also a basis for $\mathcal{O}_F$. Eventually, one has
		\begin{align*}
			G_{TB'}=B'^TT^TTB'=U^TB^T\diag{\delta^2,\sigma(\delta^2),\sigma^2(\delta^2)}BU=G_{B_{\mathfrak I'}}.
		\end{align*} By assumption,  the matrix $G_{B_{\mathfrak I'}}$ satifies the conditions in   \Cref{thm:wr3_nec_suff}. Thus,  the lattice obtained after twisting $B'$ by $T$ is WR, and it is similar to the ideal lattice $\mathfrak I$.
	\end{proof}

	\begin{lemma}\label{lemma:unit_u'}
		Let $F$ be a cyclic cubic field with a unit element that is neither totally positive nor totally negative. Assume that  $\mathcal{O}_F$ has a totally positive unit  $u$, then there exists a unit $u'$ such that $u = u'^2$.
	\end{lemma}
	\begin{proof}		
        It is known that $\mathcal{O}_K$ has fundamental units of the form $\varepsilon_1, \varepsilon_2=\sigma(\varepsilon_1)$ for some unit $\varepsilon_1$ by \cite{hasse50}. If $\varepsilon_1$ is totally positive or totally negative, then so is $\varepsilon_2$, hence $F$ cannot have a unit that is neither totally positive nor totally negative, as in the assumption. Thus $(\varepsilon_1, \sigma(\varepsilon_1), \sigma^2(\varepsilon_1)) \in \mathbb{R}^3$ must have one of the sign $\pm (+, +, -), \pm (+, -, -), \pm (+, -, +)$. By replacing $\varepsilon_1$ with $-\varepsilon_1$ or with $\pm \sigma(\varepsilon_1)$ or with $\pm \sigma^2(\varepsilon_1)$ if necessary, we obtain fundamental units $\varepsilon _1, \varepsilon_2$ such that  $\varepsilon _1>0$ and  $\varepsilon_2=\sigma(\varepsilon_1)<0$.\\

Since $u$ is totally positive, there exist some integers $k, \ell$ such that
	\begin{align*}
			u &= \varepsilon_1^k\varepsilon_2^\ell>0,\\
			\sigma(u) &= (\sigma(\varepsilon_1))^k (\sigma(\varepsilon_2))^\ell>0,\\
			\sigma^2(u) &= (\sigma^2(\varepsilon_1))^k (\sigma^2(\varepsilon_2))^\ell > 0.
		\end{align*}
		From \(\varepsilon_1 > 0,\varepsilon_2 < 0\), we infer from the first inequality that \(\ell\) is even. Since \(\ell\) is even, we infer from the second inequality that \(n\) is even. Hence, there exists a unit $u'$ such that $u = u'^2$.\\
	\end{proof}
	
The converse of the \Cref{prop:WRtwist_from_prin_ideal} holds under the following additional assumptions..
	
\begin{theorem}\label{thm:prin_ideal_from_WRtwist}
	Let $F$ be a cyclic cubic field with a unit element that is neither totally positive nor totally negative. Suppose that $\mathcal{O}_F$ has a WR twist lattice $L$ with a good basis $B=\{x,y,z\}$. Let $\alpha_0$ as in \eqref{eq:alpha2_cylic}. If $\alpha_0 = k\psi$ with $\Nm(\psi)\mid \Delta_F^t$ for some $t,k\in \Z_{>0}$ and $\psi\in \mathcal{O}_F$, then the principal ideal generated by $\psi^2$ is similar to $L$. Moreover, the Gram matrix of every minimal basis of the WR twist $L$ of $\mathcal{O}_F$ is of the form
        \begin{align*}
        \begin{bmatrix}
			s & u & v \\
			u & s & w \\
			v & w & s
		\end{bmatrix}
    \end{align*}
    where $|u|=|v|=|w|.$
\end{theorem}
\begin{proof}
Let $\alpha_0,\beta_0,\gamma_0$ as in \eqref{eq:alpha2_cylic}. Then these elements belong to $\mathcal{O}_F$ and $\beta_0=\sigma(\alpha_0), \gamma_0 = \sigma(\beta_0)$. Hence, they all have the same sign by the assumption that $B$ is a good basis and by \Cref{rem:twist_ideal_cyclic}. Multiplying by $-1$ if necessary, we can assume that they are all positive. By \Cref{lem:squareideal}, we obtain $(\alpha_0)=k(\psi)=kq\circbrack{\psi^4}$ for some $q \in \Q_{>0}$. It implies that there exists a unit element $\varepsilon'$ such that $\alpha_0= kq\varepsilon'\psi^4$. Thus $\varepsilon'$ is totally positive.  By  \Cref{lemma:unit_u'}, $\varepsilon'=\varepsilon^2$ for some unit element $\varepsilon$.  Let $\mathfrak J=\circbrack{\psi^2}$. Then $B_\mathfrak J=\set{\varepsilon\psi^2 x, \varepsilon\psi^2y,\varepsilon\psi^2y}$ is a basis of $\mathfrak J$.

 Moreover, one has \begin{align*}
G_{B_\mathfrak J}&=B_\mathfrak J^T B_\mathfrak J= B^T\text{diag}\circbrack{\varepsilon^2\psi^4,\sigma\circbrack{\varepsilon^2\psi^4},\sigma^2\circbrack{\varepsilon^2\psi^4}}B\\&= \dfrac{1}{kq}B^T\text{diag}\circbrack{\sqrt{\alpha_0},\sqrt{\beta_0},\sqrt{\gamma_0}}B. 
\end{align*}
Therefore, the ideal lattice $\mathfrak J$ is similar to $L$.

    Since $\Nm\circbrack{\psi}\mid \Delta_F^t$, the ideal factorization of $(\psi)$ is contains only ramified prime ideals. Thus, $\mathfrak J=\circbrack{\psi^2}= l\cdot \mathfrak I$ for some $l\in \Z^+$ and integral ideal $\mathfrak I$ such that $\Nm(\mathfrak I)\mid \Delta_F$. It follows that $\mathfrak I$ is similar to $L$. The remaining part of the proposition then follows from \Cref{cor:Gram_matrix}. 
	\end{proof}

\begin{corollary}
Let $F$ be a cyclic cubic field with a unit element $\epsilon$ that is neither totally positive nor totally negative, $B=\set{x,y,z}$ be a good basis of $\mathcal{O}_F$ with $\alpha_0,\beta_0,\gamma_0$ as in \eqref{eq:initial_alphabeta}. Suppose that $\ell$ is the largest positive integer such that $\dfrac{\alpha_0}{\ell}\in \mathcal{O}_F$. Let $T=\diag{\sqrt{|\alpha_0|},\sqrt{|\beta_0|},\sqrt{|\gamma_0|}}$ and let  $L$ be the twist lattice obtained by twisting this basis $B$ by $T$. Then the basis $TB$ of $L$ has the Gram matrix of the form $\ell A$, where $A$ is a matrix of the form given in \Cref{thm:wr3_nec_suff}. \\
Moreover, if $\det A\mid \Delta_F^t$ for some $t\in \Z_{>0}$, then there exists an integral ideal $\mathfrak I$ such that the ideal lattice of $\mathfrak I$ is similar to $L$. 

\end{corollary}
	
Finally, we provide some conditions for $\mathcal{O}_F$ to have an orthogonal WR twist.

\begin{theorem}\label{thm:ortho_twist}
Let $K$ be a totally real number field of degree $n\leq 7$ over $\Q$.  Suppose that the different ideal $\frak{d}_K$ of $K$ is principal.  If there exists a totally positive generator $\delta$ of $\frak{d}_K$, then the orthogonal lattice $\Z^n$ is a twist of its ring of integers $\mathcal{O}_K$.
\end{theorem}
	
\begin{proof}
Let $\Delta_K$ be the discriminant of $K$, which is the determinant of the trace form $(x,y)\mapsto \mathrm{Tr}(xy)$ on $\mathcal{O}_K$.  Moreover, we have $\Delta_K = \Nm_{K/\Q}(\delta)$.  By definition of the different ideal, we have that the twisted trace form $(x,y)\mapsto \mathrm{Tr}(\delta^{-1}xy)$ is integral and has determinant $\Delta_K/\Nm_{K/\Q}(\delta) = 1$.  
    Moreover, if $\delta_i$ denotes the image of $\delta$ under the $i$th embedding of $K$ into $\R$, then all $\delta_i > 0$ and this trace form is the inner product on the twisted lattice $T_{\delta^{-1}}\mathcal{O}_K = \mathrm{diag}\left(\sqrt{\delta_i^{-1}}\right)\cdot \mathcal{O}_K$.  The lattice $T_{\delta^{-1}}\mathcal{O}_K$ is a (positive definite, integral) unimodular lattice, and in dimension $n\leq 7$ the only such lattice is $\Z^n$ (see \cite[Section 2.4, Table 2.2]{conwaysloane1999}).
\end{proof}
	
In case of cyclic cubic fields, $\delta$ as in \Cref{thm:ortho_twist} is not necessarily positive. It is also a corollary of \Cref{prop:WRtwist_from_prin_ideal} 
    \begin{proposition}\label{prop:dif_ideal_pri}
       If the different ideal of a cyclic cubic field is principal, then its ring of integers has an orthogonal WR twist.
    \end{proposition}
\begin{proof}
Let $F$ be a cyclic cubic field with the conductor $m$ and the principal different ideal $\mathfrak{d}_F$. If $3\nmid m$, then the ideal $\mathfrak{d}_F$ coincides with the orthogonal and WR ideal as in \cite[Proposition 3.7]{tran2023well}. If $3\mid m$, then $\mathfrak{d}_F$ is similar to the orthogonal and WR ideal given in \cite[Proposition 3.18]{tran2023well}. The result is implied by using \Cref{prop:WRtwist_from_prin_ideal}.
\end{proof}

The below is a consequence of \Cref{prop:dif_ideal_pri} and \Cref{lem:dif_idealPID}. 
\begin{corollary}\label{prop:cyc_cu_orth}
    Let $F$ be a cyclic cubic field of conductor $m$, defined as in \Cref{sec:cyclic_cubic_fields}. If $b=3$ or $m$ is prime, then $\mathcal{O}_F$ has an orthogonal WR twist lattice.	
\end{corollary}

    %%%%%%%%%%%%%%%%%%%%%%%%%%%%%%%%%%%%%
\section{Explicit WR twists of the rings of integers of several families of cyclic cubic fields} \label{sec:explicit_wr_twists}
    
	In this section, we present several WR twists of the rings of integers from three families of cyclic cubic fields by Shanks ~\cite{shanks1974simplest}, Washington \cite{washington1997family}, and Kishi \cite{kishi2003family}. 
    We show that Shanks's family \cite{shanks1974simplest} has infinitely many fields of which the ring of integers has an orthogonal WR twist (\Cref{prop:infinite_ortho}).  
    The other results in Propositions \ref{prop:wr_non_ortho_shank}, \ref{prop:wash_good_bases}, and \ref{prop:kishi_good_bases} are obtained by applying  \Cref{rem:twist_ideal_cyclic}. 
    %    Even though the four steps discussed in \Cref{rem:twist_ideal_cyclic} do not form a complete algorithm for computing all good bases, they can be used as a practical method to search for some of them. 
    Starting from a known integral basis (shown in Lemmata \ref{lem:wash_disc_neven}, \ref{lem:wash_disc_nodd}, and \ref{lem:kishi_integral_basis}), we then generate new ones by applying some random unimodular transformations. After that, we test whether each resulting basis is good or not employing \Cref{rem:twist_ideal_cyclic}.

The good bases obtained using \Cref{rem:twist_ideal_cyclic} for some cyclic cubic fields from Shanks, Washington, and Kishi's families are summarized in Propositions \ref{prop:wr_non_ortho_shank}, \ref{prop:wash_good_bases}, and \ref{prop:kishi_good_bases}, respectively. Each proposition is verified following the outlined procedure: we use \texttt{SageMath} to compute and confirm the unimodularity of the transition matrices, apply steps 1–3 of \Cref{rem:twist_ideal_cyclic} (see \href{https://github.com/hoainam-le/wril_cubic_fields/blob/main/notebooks/sagemath/compute_gram_matrix_of_good_bases.ipynb}{here}), and use \texttt{Mathematica} for steps 2 and 4 of \Cref{rem:twist_ideal_cyclic} (\href{https://github.com/hoainam-le/wril_cubic_fields/blob/main/output/wolfram/solve_ineqs_same_sign_conditions.pdf}{here} and \href{https://github.com/hoainam-le/wril_cubic_fields/blob/main/output/wolfram/solve_ineqs_wr_conditions.pdf}{here}). 
Therefore, for brevity, we only present the Gram matrices of the good bases after twisting; detailed calculations are provided in the GitHub repository \cite{github_computations_wr}. In addition, for every good basis from these families, we also check whether there exists a principal ideal similar to the corresponding twist of that basis, and explicitly determine the ideal if it exists. The results are presented in \Cref{cor:wash_principal_link_with_twist}.

We remark that cyclic cubic fields have the following property.
\begin{lemma}\label{lem:totally_pos_generat}
	Let \( F \) be a cyclic cubic field, and suppose that \( F \) contains a unit that is neither totally positive nor totally negative. Then every principal ideal of \( F \) is generated by a totally positive element.
\end{lemma}
\begin{proof}
	Suppose \(\mathfrak I \) is a principal ideal of \( F \) generated by \( \delta \). Let \( u \) be a unit that is neither totally positive nor totally negative. If \( \delta / u \) is not totally positive, then we can replace \( u \) with \( -u \) or with \( \pm \sigma(u) \) or \( \pm \sigma^2(u) \) so that \( \delta / u \) is totally positive. Hence, \(\mathfrak I \) is generated by a totally positive element.
\end{proof}

    \begin{remark}\label{rmk:totally_pos_generat}
       Cyclic cubic fields from all three families Shanks ~\cite{shanks1974simplest}, Washington \cite{washington1997family}, and Kishi \cite{kishi2003family} contain units that are neither totally positive nor totally negative. This can be briefly explained by noting that the three roots of the polynomials defining these families are all units and are conjugate. Furthermore, the elementary symmetric polynomials in these roots can be determined from the defining polynomial, and \Cref{lem:sign_of_three_numbers} can be applied to verify that they do not all have the same sign. 
       Therefore, for fields in these three families, we can conclude by \Cref{lem:totally_pos_generat} that any of their principal ideals is generated by a totally positive element. However, for general cyclic cubic fields, there exist examples in which every unit is either totally positive or totally negative. Such examples were found using \texttt{SageMath} (see \href{https://github.com/hoainam-le/wril_cubic_fields/blob/main/notebooks/sagemath/all_fund_units_are_totally_positive.ipynb}{here}).

    \end{remark}

%%%%%%%%%%%%%	
\subsection{Shanks's simplest cubic fields}\label{sec:shank_family}
	
	Consider the simplest cyclic cubic field $F$ by  Shanks \cite{shanks1974simplest} defined by a polynomial \begin{align}\label{eq:shank_dfpol}
		df(x) = x^3 -nx^2-(n+3)x-1\  \text{for}\ n\in \Z,n\ge -1.
	\end{align}
	
Let $\rho$ be a root of $df(x)$, and set $F = \mathbb{Q}(\rho)$. 
Let $\mathfrak{C}_F$ be the unique integral ideal of $\mathcal{O}_F$ whose norm equals the conductor of $F$. 
By \cite[Theorem~1.4]{kashio21}, we have necessary and sufficient conditions for $\mathfrak{C}_F$ to be principal. 
Note that the different ideal of $F$ satisfies $\mathfrak{d}_F = \mathfrak{C}_F^2$. Combining \cite[Theorem~1.4]{kashio21} with \Cref{lem:I_square_prin}, 
we obtain necessary and sufficient conditions for the different ideal $\mathfrak{d}_F$ to be principal. 
When $\mathfrak{d}_F$ is principal, we may either apply \Cref{thm:ortho_twist} and \Cref{rmk:totally_pos_generat}, or  \Cref{prop:dif_ideal_pri}, to conclude that $\mathcal{O}_F$ admits an orthogonal WR twist lattice.

\begin{proposition}
    \label{prop:infinite_ortho}
    Let $F$ be a cyclic cubic field from Shanks's family as defined in \eqref{eq:shank_dfpol}. If $n^2 +3n +9$ is squarefree, then $\mathcal{O}_F$ has an orthogonal WR twist lattice. In particular, there are infinitely many simplest cubic fields of which the ring of integers has an orthogonal WR twist lattice.
   \end{proposition}
\begin{proof}
   For a simplest cubic field $F$, by \cite[Lemma 1]{cusick84}, if we choose $n$ such that $n^2 +3n +9$ is squarefree, then the discriminant of the defining polynomial in \eqref{eq:shank_dfpol} is equal to the discriminant of the field. It follow that $\mathcal O_F = \mathbb Z[\rho]$. By  \cite[Theorem 4.3]{conrad2009different}, the different ideal of $F$ is principal. It follows that $\mathcal{O}_F$ has an orthogonal WR twist by \Cref{prop:dif_ideal_pri}.
   
   The last statement in the proposition follows from the fact that there are infinitely many such integers $n$ such that $n^2 +3n +9$ is squarefree \cite[Lemma 2]{cusick84}.
\end{proof}

On the other hand, if the different ideal is not principal, then by \cite[Theorem 4.3]{conrad2009different}, $\mathcal{O}_F$ is not monogenic. In this case, \cite[Proposition 3.1]{gil2025additive} states that $n \equiv 3, 21 \pmod{27}$, $n > 12$, and $\frac{\Delta_F}{27}$ is squarefree if and only if  $\mathcal{O}_F$ has an integral basis of the form $\left\{1, \rho, \dfrac{1 + \rho + \rho^2}{3}\right\}$. Using this basis, we can explicitly construct a good basis as below.
	
\begin{proposition}\label{prop:wr_non_ortho_shank}
	Let $F$ be the simplest cubic field defined by a polynomial as in \eqref{eq:shank_dfpol} with an integral basis of the form $\left\{1, \rho, \dfrac{1 + \rho + \rho^2}{3}\right\}$. Then $\left\{\dfrac{1+\rho+\rho^2}{3},\rho,\rho+\rho^2\right\}$ is a good basis of $\mathcal{O}_F$.
\end{proposition}
	
\begin{proof}Here is the Gram matrix of this basis after twisting.
		\begin{align*}
			\|Tx\|^2 &= \dfrac{1}{9} (n^2+3n+3)(n^2+3n+9),\\
			\langle Tx,Ty\rangle  &= \dfrac{-1}{27}(n^2+9n+9)(n^2+3n+9),\\
			\langle Tx,Tz\rangle &= \dfrac{1}{27}(n^2-3n-9)(n^2+3n+9),\\
                \langle Ty,Tz\rangle &= \dfrac{2}{3}(n^2+3n+9).
		\end{align*}
	\end{proof}
	
\begin{remark}
	\Cref{prop:wr_non_ortho_shank} presents cases in which $\mathcal{O}_F$ yields a non-orthogonal WR lattice. Moreover, the orthogonal WR twist given in \Cref{thm:ortho_twist} is not unique. This means that even if the different ideal is principal, $\mathcal{O}_F$ may still admit another non-orthogonal WR twist. For example, consider the case $n = 39$ (see \href{https://github.com/hoainam-le/wril_cubic_fields/blob/main/notebooks/sagemath/shanks_twist_example_n39.ipynb}{here}).
        In this case, we have $[\mathcal{O}_F : \Z[\rho]] = 27$, and by \cite[Corollary 1.6]{kashio21}, the field~$F$ is monogenic. It follows that the different ideal $\frak{d}_F$ is principal. By \Cref{prop:dif_ideal_pri}, $\mathcal{O}_F$ therefore has an orthogonal WR twist. Furthermore, $\mathcal{O}_F$ has an integral basis
	\begin{align*}
			\left\{\dfrac{-2\rho^2+\rho+1}{9},\dfrac{4\rho^2+4\rho+1}{9},\dfrac{-\rho^2-2\rho}{3}\right\},
	\end{align*}
		whose Gram matrix satisfies the conditions of \Cref{thm:wr3_nec_suff} and gives rise to a non-orthogonal WR twist.
	\end{remark}

%%%%%%%%%%%%%%%%%%%%
\subsection{Washington's cyclic cubic fields} \label{sec:washington_family}
	
Consider a cyclic cubic field $F$, by Washington  \cite{washington1997family}, which defined by a polynomial \begin{align}\label{eq:wash_df}
	    df(x) =  x^3 -\circbrack{n^3-2n^2+3n-3}x^2-n^2x-1.
	\end{align}
	For all $n\ne 1$, the polynomials $df(x)$ is irreducible with the discriminant \begin{align*}
		\Delta_{df}=\circbrack{n-1}^2\circbrack{n^2+3}^2\circbrack{n^2-3n+3}^2.
	\end{align*}
	Let $\rho$ be a root of $df(x)$. The Galois group is generated by the transformation 
    \begin{align*}
		\rho \mapsto -\dfrac{\rho+1}{\circbrack{n^2-n+1}\rho+n}.
	\end{align*}
	According to the comment at the beginning of \cite[Section 2]{washington1997family}, $\circbrack{n-1},\circbrack{n^2+3},\circbrack{n^2-3n+3}$ are three numbers that are pairwise relatively prime, except possibly for powers of 2 and 3. From \cite[Theorem 1]{washington1997family}, we write 
    \begin{align}\label{eq:disc}
		\circbrack{n^2+3}\circbrack{n^2-3n+3}=bc^3
	\end{align} 
    where $b$ is cubefree, then the discriminant 
	\newcommand{\disc}{\text{disc}}
	\begin{align*}
		\Delta_F= 81^\delta \prod_{p\mid b, p\ge 5}p^2, 
	\end{align*}where $\delta =0$ if $n\not\equiv 0 \pmod 3$ and $\delta =1$ if $n\equiv 0\pmod 3.$
	
In $F$, some elements listed in the lemma below are algebraic integers. This can be easily verified using \texttt{SageMath} (see \href{https://github.com/hoainam-le/wril_cubic_fields/blob/main/notebooks/sagemath/algebraic_elements.ipynb}{here}) by computing their minimal polynomials. These elements will later be shown to form parts of an integral basis of $F$.
\begin{lemma}\label{lem:wash_integer}
	For all $n$, the element $\dfrac{\rho^2-1}{n-1}$ is in $\mathcal{O}_F$. Moreover, if $n$ odd, then $\dfrac{\rho^2-1}{2n-2}$ and $\dfrac{\rho^2+\rho}{2}$ are in $\mathcal{O}_F$.
\end{lemma}

\begin{lemma}\label{lem:wash_disc_neven}
	If $n$ is even and $\dfrac{\circbrack{n^2+3}\circbrack{n^2-3n+3}}{9^\delta}$ is squarefree, then $\Delta_F=\circbrack{n^2+3}^2\circbrack{n^2-3n+3}^2$. Moreover, $\set{1,\rho,\dfrac{\rho^2-1}{n-1}}$ is an integral basis of $\mathcal O_{F}$.
\end{lemma}
\begin{proof} If $\dfrac{\circbrack{n^2+3}\circbrack{n^2-3n+3}}{9^\delta}$ is squarefree, then $c$ in \eqref{eq:disc} must be $1$. It implies that 
    $$\circbrack{n^2+3}\circbrack{n^2-3n+3}=b,\ \text{where}\ b\ \text{is cubefree}.$$ 
		By \cite[Theorem 1]{washington1997family}, it's sufficient to prove taht $\dfrac{\circbrack{n^2+3}\circbrack{n^2-3n+3}}{9^\delta}$ is not divisible by $3$ and $2$. Since $n$ is even,  $b=\circbrack{n^2+3}\circbrack{n^2-3n+3}$ (as in \eqref{eq:disc}) is odd.
		
		\begin{itemize}
			\item If $3\nmid n$, then $\delta =0$ and $3\nmid b$. It follows that $b^2 = \prod_{p\mid b, p\ge 5}p^2=\Delta_F$.
			\item If $3\mid n$ and we write $n=3k$, then $b=  9\circbrack{3k^2+1}\circbrack{3k^2-3k+1}$. Moreover, it is clearly that $3\nmid \circbrack{3k^2+1}\circbrack{3k^2-3k+1}$. It implies that $b=  9\prod_{p\mid b, p\ge 5}p$ and therefore $$b^2 = 81\prod_{p\mid b, p\ge 5}p^2=\Delta_F.$$
		\end{itemize} 
        Additionally, we have\begin{align*}
			\begin{vmatrix}
				1&\rho& \dfrac{\rho^2-1}{n-1}\\
				1&\sigma(\rho)&\dfrac{\sigma(\rho^2)-1}{n-1}\\1&\sigma^2(\rho)&\dfrac{\sigma^2(\rho^2)-1}{n-1}
			\end{vmatrix}^2=  \dfrac{1}{\circbrack{n-1}^2}\begin{vmatrix}
				1&\rho&\rho^2\\
				1&\sigma(\rho)&\sigma(\rho^2)\\
				1&\sigma^2(\rho)&\sigma^2(\rho^2)
			\end{vmatrix}^2&=\dfrac{1}{(n-1)^2}\Delta_{df}\\
			&=\circbrack{n^2+3}^2\circbrack{n^2-3n+3}^2= \Delta_F.
		\end{align*} Hence, $\set{1,\rho,\dfrac{\rho^2-1}{n-1}}$ is an integral basis of $\mathcal{O}_F$. 
	\end{proof}
	
	\begin{lemma}\label{lem:wash_disc_nodd}
		If $n$ is odd and $\dfrac{\circbrack{n^2+3}\circbrack{n^2-3n+3}}{4\cdot9^\delta}$ is squarefree, then $\Delta_F= \dfrac{\circbrack{n^2+3}^2\circbrack{n^2-3n+3}^2}{16}$. Furthermore, $\set{\dfrac{\rho^2-1}{2n-2}, \dfrac{\rho^2+\rho}{2},\rho^2}$ is an integral basis of $\mathcal{O}_F$.
	\end{lemma}
	\begin{proof}
		Note that in this case, one has $b=\circbrack{n^2+3}\circbrack{n^2-3n+3}$ (in \eqref{eq:disc})  which is cubefree. When $n$ is odd, write $n=2k+1$, then we obtain that  
        \begin{align*}
			b= 4\circbrack{k^2+k+1}\circbrack{4k^2-2k+1}.
		\end{align*} It is easy to prove that the last two factors of $b$  are odd. Therefore, one has $\dfrac{b}{4\cdot 9^\delta }$ is not divisible by $2$. By a similar argument of proof of  \Cref{lem:wash_disc_neven}, $\dfrac{b}{4\cdot 9^\delta }$ is also not divisible by $3$. Then we have $b=4\cdot 9^\delta \prod_{p\mid b, p
			\ge 5}p$. As a consequence, $\Delta_F=\dfrac{b^2}{16}.$    \\

        Similarly, one has
        \begin{align*}
         \begin{vmatrix}
            \dfrac{\rho^2-1}{2n-2}& \dfrac{\rho^2+\rho}{2}&\rho^2\\
            \dfrac{\sigma(\rho)^2-1}{2n-2}& \dfrac{\sigma(\rho)^2+\sigma(\rho)}{2}&\sigma(\rho)^2\\
             \dfrac{\sigma^2(\rho)^2-1}{2n-2}& \dfrac{\sigma^2(\rho)^2+\sigma^2(\rho)}{2}&\sigma^2(\rho)^2
        \end{vmatrix}^2&=\dfrac{1}{16\circbrack{n-1}^2}\begin{vmatrix}
            1&\rho&\rho^2\\
            1&\sigma(\rho)&\sigma(\rho^2)\\
            1&\sigma^2(\rho)&\sigma^2(\rho^2)
        \end{vmatrix}=\dfrac{\circbrack{n^2+3}^2\circbrack{n^2-3n+3}^2}{16}.
    \end{align*}
    Therefore, $\set{\dfrac{\rho^2-1}{2n-2}, \dfrac{\rho^2+\rho}{2},\rho^2}$ is also an integral basis of $\mathcal O_F$.
\end{proof}

 The determinant computations in Lemmata \ref{lem:wash_disc_neven} and \ref{lem:wash_disc_nodd} can also be done by using  \texttt{SageMath} (see \href{https://github.com/hoainam-le/wril_cubic_fields/blob/main/notebooks/sagemath/compute_disc.ipynb}{here}). Moreover, from these lemmata, we observe that if $(n^2-3n+3)$ divides the discriminant, then there is a unique ideal of norm $(n^2-3n+3)$. Furthermore,  we have 
 $$\Nm\circbrack{\dfrac{\rho^2-1}{n-1}-n\rho} = -(n^2-3n+3).$$
 Thus, one obtains the following lemma.

\begin{lemma}\label{lem:princ_gen}
If \( n^2 - 3n + 3 \) divides the discriminant of \( F \), then the unique ideal of norm  \( n^2 - 3n + 3 \) in \( \mathcal{O}_F \) is principal and generated by
\(
\frac{\rho^2 - 1}{n - 1} - n\rho.
\)
\end{lemma}
Explicit good bases of $\mathcal{O}_F$ are given in the following proposition.
\begin{proposition}\label{prop:wash_good_bases}
    \begin{enumerate}
        \item If $n$ is as in \Cref{lem:wash_disc_neven} and $n\ne 2$, then $\set{\rho, \dfrac{\rho^2-1}{n-1}-\rho, \rho^2}$ is a good basis of $\mathcal{O}_F$.
        \item  If $n$ is as in \Cref{lem:wash_disc_nodd}, then
        \begin{enumerate}[label = \alph*.]
                \item $\set{\rho^2,\dfrac{n-2}{2n-2}\rho^2+\dfrac{1}{2}\rho+\dfrac{1}{2n-2},\dfrac{\rho^2+\rho}{2}}$ is a good basis of $\mathcal O_F$ when $n\ge 5$.
                \item $\set{\dfrac{n-2}{2n-2}\rho^2+\dfrac{1}{2}\rho+\dfrac{1}{2n-2},-\dfrac{\rho^2-1}{2n-2}+\rho,\dfrac{\rho^2+\rho}{2}}$ is a good basis of $\mathcal O_F$ when $n\ge 0$.
            \end{enumerate}
    \end{enumerate}
\end{proposition}

    \begin{proof}
We only show the Gram matrices corresponding to good bases of $\mathcal{O}_F$. See the computation using \texttt{SageMath} at \href{https://github.com/hoainam-le/wril_cubic_fields/blob/main/notebooks/sagemath/compute_gram_matrix_of_good_bases.ipynb}{here}.
\begin{enumerate}
    \item \begin{align*}
         \|Tx\|^2 &= (n - 1)^2  (n^2 - 3n + 3)  (n^2 - n + 3) (n^2 + 3),\\
\langle Tx,Ty\rangle  &= n (n - 1)^2  (n^2 - 3n + 3)  (n^2 + 3),\\
\langle Tx,Tz\rangle &= - n  (n - 1)^2  (n^2 - 3n + 3)  (n^2 + 3),\\
\langle Ty,Tz\rangle &=  n  (n - 1)^2  (n^2 - 3n + 3)  (n^2 + 3).&
\end{align*}
    \item 
    \begin{enumerate}[label = \alph*.]
        \item \begin{align*}
         \|Tx\|^2 &= \dfrac{1}{16} (n^2-3n+3)(n^2+3)(n^4 - 5n^3 + 10n^2 - 11n + 1),\\
\langle Tx,Ty\rangle  &= \dfrac{1}{32}(n^2 - 3n + 3)( n^2 + 3)(n^2 - 2n - 1)(n^2 - 4n + 7),\\
\langle Tx,Tz\rangle &= \dfrac{1}{32}( n^2 - 3n + 3)(n^2 + 3)(n^4 - 8n^3 + 16n^2 - 16n - 1),\\
\langle Ty,Tz\rangle &= \dfrac{1}{64}(n^2 - 3n + 3 )( n^2 + 3)( n - 1)(n^3 - 11n^2 + 19n - 1).&
\end{align*}
\item 
\begin{align*}
 \|Tx\|^2 &= \dfrac{1}{32} (n^2 - 4n + 7)(n^2 - 2n + 3)(n^2 - 3n + 3)( n^2 + 3 ),\\
\langle Tx,Ty\rangle  &= \dfrac{1}{64}(n - 3 )( n - 1)(n^2 - 4n + 7)(n^2 - 3n + 3)(n^2 + 3),\\
\langle Tx,Tz\rangle &= \dfrac{1}{64}(n - 3 )( n - 1)(n^2 - 4n + 7)(n^2 - 3n + 3)(n^2 + 3),\\
\langle Ty,Tz\rangle &= \dfrac{1}{64}(n - 3 )( n - 1)(n^2 - 4n + 7)(n^2 - 3n + 3)(n^2 + 3).&
\end{align*}
    \end{enumerate}
\end{enumerate}
    \end{proof}

From \Cref{prop:wash_good_bases}, we observe that the Gram matrices in Cases 1 and 2b match the form described in \Cref{thm:prin_ideal_from_WRtwist}. Particularly, they can be rewritten as follows:
    \begin{align*}
    (n-1)^2(n^2+3)(n^2-3n+3)
    \begin{bmatrix}
    n^2 - n + 3 & n & -n \\
    n & n^2 - n + 3 & n \\
    -n & n & n^2 - n + 3
    \end{bmatrix},
    \end{align*}
    and
    \begin{align*}
    \dfrac{1}{64}(n^2 - 4n + 7)(n^2 + 3)(n^2 - 3n + 3)
    \begin{bmatrix}
    2n^2 - 4n + 6 & (n - 3)(n - 1) & (n - 3)(n - 1) \\
    (n - 3)(n - 1) & 2n^2 - 4n + 6 & (n - 3)(n - 1) \\
    (n - 3)(n - 1) & (n - 3)(n - 1) & 2n^2 - 4n + 6
    \end{bmatrix},
    \end{align*}
respectively.

Furthermore,  for each good basis, we can compute $\alpha^2$ using the formula \eqref{eq:alpha2_cylic} by \texttt{SageMath} (see \href{https://github.com/hoainam-le/wril_cubic_fields/blob/main/notebooks/sagemath/compute_gram_matrix_of_good_bases.ipynb}{here}), and then obtain the following results:
\begin{itemize}
\item For the first good basis in Case 1, we have
\begin{align*}
    \alpha^2=\circbrack{n-1}^2\circbrack{n^2+3}\psi \text{ where }\text{N}\circbrack{\psi}=\circbrack{n^2-3n+3}^2\mid \Delta_F.
\end{align*}
\item For the second good basis  in Case 2b, we have
\begin{align*}
    \alpha^2=\circbrack{n^2-4n+7}\circbrack{n^2+3}\psi\text{ where }\text{N}\circbrack{\psi}=\circbrack{n^2-3n+3}^2\mid \Delta_F.
\end{align*}
\end{itemize}

Therefore, by combining this with \Cref{lem:princ_gen} and \Cref{thm:prin_ideal_from_WRtwist}, the corollary is implied.

\begin{corollary} \label{cor:wash_principal_link_with_twist}
    \begin{enumerate}
    \item If $n$ is as in \Cref{lem:wash_disc_neven} and $n\ne 2$, then the lattice obtained by twisting the good basis in \Cref{prop:wash_good_bases}, Case 1 is similar to the unique WR principal ideal of norm $n^2-3n+3$.
    
    \item If $n$ is as in \Cref{lem:wash_disc_nodd} and $n\ge 0$, then the lattice obtained by twisting the good basis in \Cref{prop:wash_good_bases}, Case 2b is similar to the unique WR principal ideal of norm $n^2-3n+3$.
    \end{enumerate}
\end{corollary} 
\begin{proof} Recall that $\mathcal{O}_F$ always contains a unit element that is neither totally positive nor totally negative. Let $m$ be the conductor of $F$. Then $\Delta_F=m^2$.
    \begin{enumerate}
    \item In this case, since $(n^2 - 3n + 3)$ is a squarefree divisor of $m$, it follows from \cite[Theorem 1.2]{tran2023well} that the unique  ideal of norm $n^2 - 3n + 3$ is WR if and only if $\dfrac{n^2+3}{4}\le n^2-3n+3\le 4(n^2+3)$. The inequalities are satisfied for $n \ne 2$. Applying \Cref{thm:prin_ideal_from_WRtwist}, the result follows immediately.
    
    \item Similarly, the inequalities $\dfrac{n^2+3}{16}\le n^2-3n+3\le n^2+3$ holds for $n\ge 0$. The result is then followed from \cite[Theorem 1.2]{tran2023well}.
    \end{enumerate}
\end{proof}

%%%%%%%%%%%%%%%%    
\subsection{Kishi's Cyclic cubic fields} \label{sec:kishi_family}
The family of cyclic cubic fields by Kishi \cite{kishi2003family} is defined by  polynomials of the form \begin{align}
		df(x)=x^3-n(n^2+n+3)(n^2+2)x^2-(n^3+2n^2+3n+3)x-1,n\in \Z
	\end{align}
	with discriminant
	\begin{align*}
		\Delta_{df} = (n^2+1)^2(n^2+3)^2(n^4+n^3+4n^2+3)^2.
	\end{align*}
Let $\rho$ be a root of $df(x)$. Using the table at the beginning of \cite[Section 4]{balady2016families}, one has the Galois group generated by the transformation. \begin{align*}
		\rho \mapsto -\dfrac{n\rho+1}{\circbrack{n^4 + n^3 + 3n^2 + n + 1}\rho+(n^2 + n + 1)}.
	\end{align*}
    
	Denote \begin{align}
		\label{eq:yashuiro_bc}\circbrack{n^2+3}\circbrack{n^4+n^3+4n^2+3}=ab^3
	\end{align} with $a$ is cube-free. Then the discriminant of $F$ is \begin{align*}
		\Delta_F=  \circbrack{3^\delta \prod_{p\mid a, p\ge 7}p}^2 \ \text{where }\delta =  \left\{\begin{matrix}
			0&\text{if $n\equiv 2\pmod 3$ or $n\equiv 16 \pmod {27}$}\\2&\text{otherwise}
		\end{matrix}\right..
	\end{align*}
	
	Let 
    \begin{equation}\label{eq:N}
        N =  \dfrac{\circbrack{n^2+3}\circbrack{n^4+n^3+4n^2+3}}{4^{\delta_1}9^{\delta_2}}
    \end{equation}
    where $\delta_1= \left\{\begin{matrix}
		0&\text{if $n$ is even}\\1&\text{if $n$ is odd}
	\end{matrix}\right.$ and $\delta_2 =\left\{\begin{matrix}
		0&\text{ if $n\equiv 2\pmod 3$}\\1&\text{if $n\not\equiv 2\pmod 3$}
	\end{matrix}\right.$.

  \begin{remark}\label{rmk:kishi_disc}
      From \cite[Corrolary 1.4]{kishi2003family}, if $N$ is squarefree, then the discriminant $\Delta_F$ is
    \begin{align*}
        \Delta_F = \circbrack{n^2+3}^3\circbrack{n^4+n^3+4n^2+3}^2/e^2,
    \end{align*}
where 
\begin{align*}
   e = 
\begin{cases}
1 & \text{if } n \equiv 0,2 \pmod{6} \text{ or } n \equiv 4,10 \pmod{18},\\
3 & \text{if } n \equiv 34,52 \pmod{54},\\
4 & \text{if } n \equiv 3,5 \pmod{6} \text{ or } n \equiv 1,13 \pmod{18},\\
12 & \text{if } n \equiv 7,25 \pmod{54},\\
27 & \text{if } n \equiv 16 \pmod{54},\\
108 & \text{if } n \equiv 43 \pmod{54}.
\end{cases}
\end{align*}
  \end{remark}  
We now describe the integral basis of $\mathcal{O}_F$ in each case.
\begin{lemma}\label{lem:kishi_integral_basis} Let $N$ be the number defined as \eqref{eq:N} and let
    \[
        \theta := \frac{(3n^2+n+3)\rho^2+(n^2+n+2)\rho+1}{n^2+1}.
     \]
         If $N$ is squarefree, then the ring of integers $\mathcal{O}_F$ has an integral basis $B$ given by the following table.
        \begin{center}
        \renewcommand{\arraystretch}{1.3}
            \begin{tabular}{ ||>{\centering\arraybackslash}p{3.5cm}|>{\centering\arraybackslash}p{3.5cm}|| }
            \hline\hline
$n$ & An integral basis $B$ \\
\hline
$n\equiv 0,2\pmod{6}$ & \multirow{2}{*}{$\left\{ \theta,\rho,\rho^2 \right\}$} \\
\cline{1-1}
$n\equiv 4,10 \pmod{18}$ & \\
\hline
$n\equiv 1,13\pmod{18}$ & \multirow{2}{*}{$\left\{ \frac{ \theta}{2},\frac{\rho^2+\rho}{2},\rho^2 \right\}$} \\
\cline{1-1}
$n\equiv 3,5 \pmod{6}$ & \\
\hline
$n\equiv 7,25\pmod{54}$&$\left\{ \frac{ \theta}{6},\frac{\rho^2+\rho}{2},\rho^2 \right\}$\\\hline
$n\equiv 16\pmod{54}$&$\left\{ \frac{ \theta}{9},\frac{2\rho^2+\rho}{3},\rho^2 \right\}$\\\hline
$n\equiv 34,52\pmod{54}$&$\left\{ \frac{ \theta}{3},\rho,\rho^2 \right\}$\\\hline
$n\equiv 43\pmod{54}$&$\left\{ \frac{ \theta}{18},\frac{5\rho^2+\rho}{6},\rho^2 \right\}$\\\hline\hline
\end{tabular}
        \end{center}
\end{lemma}

\begin{remark}
By using \texttt{SageMath} (see \href{https://github.com/hoainam-le/wril_cubic_fields/blob/main/notebooks/sagemath/algebraic_elements.ipynb}{here}), one can easily verify that \(  \theta \in \mathcal{O}_F \). In each of the above cases, the corresponding scaling factors (such as \( 1/2 \), \( 1/3 \), \( 1/6 \), \( 1/9 \), and \( 1/18 \)), as well as the elements \( \frac{\rho^2+\rho}{2} \), \( \frac{2\rho^2+\rho}{3} \), and \( \frac{5\rho^2+\rho}{6} \), are compatible with the congruence conditions on \( n \) and ensure that these elements lie in \( \mathcal{O}_F \).
\end{remark}
    \begin{proof}
     This lemma is proved by computing the discriminant of each basis and verifying that it matches the discriminant of the field, as given in \Cref{rmk:kishi_disc}, using \texttt{Sagemath} (see \href{https://github.com/hoainam-le/wril_cubic_fields/blob/main/notebooks/sagemath/compute_disc.ipynb}{here}).
    \end{proof}

\begin{proposition}\label{prop:kishi_good_bases}
Let $N$ be the number defined as \eqref{eq:N}. If $N$ is squarefree, then $\mathcal O_F$ has some good bases $\mathcal{GB}$ as below.\\
    \begin{center}
            \renewcommand{\arraystretch}{1.5}
\scalebox{0.85}{\begin{tabular}{||>{\centering\arraybackslash}p{3.0cm}|
                >{\centering\arraybackslash}p{1.5cm}|
                >{\centering\arraybackslash}p{12cm}||}
\hline$n$ & Conditions & A good basis $\mathcal{GB}$ \\\hline\hline
$n\equiv 0,2\pmod 6$ or\break $n\equiv 4,10 \pmod{18}$&$n\ne 0$&$\set{\dfrac{n\rho^2+(n^2+n+2)\rho+1}{n^2+1}, \rho,\rho^2}$\\\hline
\multirow{3}{*}{\makecell{$n\equiv 1,13\pmod{18}$ or \\ $n\equiv 3,5\pmod 6$}}&$|n|\ge 3$&$\set{-\sign(n)\dfrac{ n\rho^2+(2n^2+n+3)\rho+1}{2(n^2+1)}, \dfrac{\rho^2+\rho}{2},\sign(n)\rho}$\\\cline{2-3}
&$|n|\ne 1$&$\set{\dfrac{(n^2 -\sign(n) n + 1)\rho^2-\sign(n)(n^2+n+2)\rho-\sign(n)\cdot 1}{2(n^2+1)}, \dfrac{\rho^2+\rho}{2},-\sign(n)\rho^2} $\\\cline{2-3}
&$n=-1$&$\set{\dfrac{(n^2 -n + 1)\rho^2-(n^2+n+2)\rho- 1}{2(n^2+1)}, \dfrac{\rho^2+\rho}{2},-\rho^2}$\\\hline
\multirow{2}{*}{$n\equiv 7,25\pmod{54}$}&$n\le -8$&$\set{\dfrac{ \theta}{6},\dfrac{\rho^2+\rho}{2},\rho^2}$\\\cline{2-3}
&$n\le -2$&$\set{\dfrac{ \theta}{6},\dfrac{-n\rho^2+(2n^2-n+1)\rho-1}{6(n^2+1)},\dfrac{-n\rho^2-(n^2+n+2)\rho-1}{3(n^2+1)}}$\\\hline
\multirow{2}{*}{$n\equiv 16 \pmod{54}$ }&{\makecell{$n\ge 7$ or\\ $n\le -6$}}&$\set{\dfrac{\rho^2-\rho}{3},\rho^2,-\dfrac{ \theta}{9}}$\\\cline{2-3}&\makecell{$n\ge 5$ or\\$n\le -2$}&$\set{\dfrac{ \theta}{9},\dfrac{n\rho^2+(4n^2+n+5)\rho+1}{9(n^2+1)},\dfrac{-n\rho^2+(n^2+n+2)\rho+1}{3(n^2+1)}}$\\\hline
$n\equiv 34,52\pmod{54}$&$|n|\ge 6$&$\set{\dfrac{n\rho^2+(n^2+n+2)\rho+1}{3(n^2+1)}, \rho,\rho^2}$\\\hline
$n\equiv 43\pmod{54}$&&$\set{\rho^2,\dfrac{-\rho^2+\rho}{6},-\dfrac{ \theta}{6}}$\\\hline
\end{tabular}}
        \end{center}
    \end{proposition}
\begin{proof}
The Gram matrices of these good bases do not have the form mentioned in \Cref{thm:prin_ideal_from_WRtwist}; therefore, we cannot apply this theorem. Instead, we directly compute WR twist following the four steps in \Cref{rem:twist_ideal_cyclic}.  Detailed computations can be found \href{https://github.com/hoainam-le/wril_cubic_fields/blob/main/notebooks/sagemath/compute_gram_matrix_of_good_bases.ipynb}{here}. 
\end{proof}
%%%%%%%%%%%%%%%%%%%%%%%%%%%%%%

%%%%%%%%%%%%%%
	
\bibliographystyle{plain}
\bibliography{bibliography}
    
\end{document}